\documentclass [11pt] {article}

\usepackage {amssymb}
\usepackage {exscale}
                                      
\title {Formal GNS Construction and States \\
        in Deformation Quantization}

\author {
        {\bf
        Martin Bordemann\thanks{Martin.Bordemann@physik.uni-freiburg.de},
        Stefan Waldmann\thanks{Stefan.Waldmann@physik.uni-freiburg.de}}
        \\[3mm]
        Fakult\"at f\"ur Physik \\
        Universit\"at Freiburg \\
        Hermann-Herder-Str. 3 \\
        79104 Freiburg i.~Br., F.~R.~G\\[3mm]}

\date {FR-THEP 96/12\\
       Revised Version \\[3mm]
       November 1997
         }

\newcommand {\bfmath} [1] {{\mathbf {#1}}}

\newcommand {\BEQ} [1] {\begin {equation} \label {#1}}
\newcommand {\EEQ} {\end {equation}}      

\newcommand {\LR} {\mathbb R (\!(\lambda)\!)}
\newcommand {\NPR} {\mathbb R \langle\!\langle\lambda^*\rangle\!\rangle}
\newcommand {\CNPR} {\mathbb R \langle\!\langle\lambda\rangle\!\rangle}

\newcommand {\LC} {\mathbb C (\!(\lambda)\!)}
\newcommand {\CNPC} {\mathbb C \langle\!\langle\lambda\rangle\!\rangle}
\newcommand {\NPC} {\mathbb C \langle\!\langle\lambda^*\rangle\!\rangle}

\newcommand {\LS} [1] {{#1} (\!(\lambda)\!)}
\newcommand {\NP} [1] {{#1} \langle\!\langle \lambda^* \rangle\!\rangle}
\newcommand {\CNP} [1] {{#1} \langle\!\langle \lambda \rangle\!\rangle}
\newcommand {\field} [1] {{\mathsf {#1}}}

\newcommand {\HC} {\mathfrak H_{\mathsf C}}
\newcommand {\KC} {\mathfrak K_{\mathsf C}}

\newcommand {\WL} {\mathcal W \otimes \Lambda}
\newcommand {\W} {\mathcal W}

\newcommand {\ad} {{\rm ad}}

\newcommand {\spec} {\mathop{\rm spec}}

\newcommand {\cc} [1] {\overline {{#1}}}
\newcommand {\supp} {{\rm supp}}
\newcommand {\SP} [2] {\left\langle {#1}\, , \, {#2} \right\rangle}
\newcommand {\norm} [1] {\left\| {#1} \right\|}

\newenvironment {proof}{\small {\sc Proof:}}{{\hspace*{\fill} $\square$}}

\newtheorem {lemma} {Lemma} 
\newtheorem {proposition} {Proposition}
\newtheorem {theorem} {Theorem}
\newtheorem {corollary} {Corollary}
\newtheorem {definition} {Definition}

\begin {document}

\maketitle

\begin {abstract}

 In this paper we develop a method of constructing Hilbert spaces
 and the representation of the formal algebra of quantum observables
 in deformation quantization which is an analog of the well-known
 GNS construction for complex $C^*$-algebras: in this approach the
 corresponding positive linear functionals (`states') take their
 values not in the field of complex numbers, but in (a suitable
 extension field of) the field of formal complex Laurent series in the
 formal parameter. By using the algebraic and topological properties
 of these fields we prove that this construction makes sense and
 show in physical examples that standard representations as the Bargmann
 and Schr\"odinger representation come out correctly, both formally
 and in a suitable convergence scheme. For certain Hamiltonian functions
 (contained in the Gel'fand ideal of the positive functional) a formal
 solution to the time-dependent Schr\"odinger equation is shown to exist.
 Moreover, we show that for
 every K\"ahler manifold equipped with the Fedosov star product of Wick
 type all the classical delta functionals are positive and give rise
 to some formal Bargmann representation of the deformed algebra.

\end {abstract}

\section {Introduction}

In the programme of deformation quantization introduced by
Bayen, Flato, Fronsdal, Lichnerowicz and Sternheimer \cite{BFFLS78}
the algebra of quantum observables is considered as an associative
local formal deformation (a so-called star product $\ast$) 
of the associative commutative algebra of smooth
complex-valued functions $C^\infty (M)$ on a given symplectic manifold $M$,
such that the
first order commutator equals $i\lambda$ times the Poisson bracket
and such that complex conjugation is an antilinear involution of the
deformed algebra. The latter is equal to $C^\infty(M)[[\lambda]]$, the space
of formal power series in the deformation parameter $\lambda$ with
coefficients in $C^\infty(M)$, and the associative noncommutative 
multiplication $\ast$ is bilinear with respect to the ring 
$\mathbb C[[\lambda]]$ of formal power series in $\lambda$ with complex 
coefficients. 
On one hand the rather difficult question of existence 
of these deformations for general symplectic manifolds has positively 
been answered (DeWilde-Lecomte 1983 \cite{DL83}; Fedosov 1985 
\cite{Fed94,Fed96}; Omori-Maeda-Yoshioka 1991 \cite{OMY91}).
Moreover their classification (up to equivalence) in terms of 
formal power series with coefficients in the second de Rham 
cohomology class of the underlying symplectic manifold has 
recently been achieved (\cite [p. 204] {NT95b}, \cite {BCG96}).

On the other hand star products have the
--at first sight rather unpleasant-- feature of lacking convergence
uniform in the formal parameter which is due to the fact that they
depend on the infinite jets of the two functions which in turn can be
made as divergent as possible by Borel's Theorem (see e.g. \cite{Whi34};
see also the article by Rubio 1984 \cite{Rub84} for the commutativity of
local associative products on $C^\infty(M)$). Moreover,
the deformed algebras do not seem to have obvious
representations in some complex separable Hilbert space which is 
unsatisfactory from the physical point of view. In the past decade, however,
several people have attacked this problem: Cahen, Gutt, and Rawnsley
\cite{CGRI}, \cite{CGRII} start from the 
finite-dimensional operator algebras of geometric quantization 
in tensor powers of a very ample regular prequantum line bundle over 
a compact K\"ahler manifold and use coherent states (see \cite{Ber74},
\cite{Raw77}) to first construct star products for the Berezin-Rawnsley
symbols (\cite{Ber74}, \cite{Raw77}) for each tensor power separately. 
In a second step 
an asymptotic expansion of these star products in the inverse tensor
power is shown to define a local star product on the manifold where
the formal parameter appears as a sort of interpolation of the inverse
tensor powers. For compact Hermitean symmetric spaces they showed that the 
subspace of representative functions has a convergent star product.
See also \cite{BBEW96}, \cite{BBEW96a} for an elementary algebraic approach 
in the particular
case of complex projective space. Pflaum has studied star products and
their convergence on suitable subspaces on cotangent bundles in his thesis
\cite{Pfl95}. 
In flat $\mathbb R^{2n}$ a non-formal analogue (`twisted products') 
of star products using integral formulas can be defined on the Schwartz 
test function space (which is thereby made into an associative 
topological complex algebra) and extended to larger spaces by 
functional analytic techniques. The usual Weyl-Moyal star product is 
recovered as a asymptotic expansion in $\hbar$ (which can be introduced 
as an dilatation parameter), see \cite {Han84,Kam86,Mai86} for more 
details. Moreover, Fedosov has shown that the deformed algebra
allows a so-called asymptotic operator representation
in a complex Hilbert space
if and only if certain integrality properties of a formal index 
are satisfied which he defines as a formal analogue of the 
Atiyah-Singer index theorem and which is an invariant of the underlying
symplectic manifold (see e. g. Fedosov's book \cite{Fed96} for a detailed
exposition).

The approach of this paper is motivated by the following consideration:
in the theory of complex
$C^*$-algebras (which form one of the main mathematical 
pillars of algebraic
quantum field theory (see Haag's book \cite{Haa93} for details)
and ---in a more algebraic context--- of Connes' noncommutative 
geometry (see Connes' book \cite {Con94})) the
representing complex Hilbert spaces for a given complex $C^*$-algebra 
$A$ are constructed by means of the so-called GNS representation 
(see e. g. \cite{BR87}): roughly speaking,
any positive complex-valued linear functional 
on $A$ (i.e. which maps positive elements of $A$
to nonnegative real numbers) gives rise to a left ideal $I$ of $A$
(the Gel'fand ideal), and $A$ is canonically represented 
on the quotient space $A$ modulo $I$. Due to the positivity
of the initial functional this quotient is equipped with a 
positive definite
sesquilinear form and thus becomes a complex pre-Hilbert space 
whose completion
yields the desired representation of the algebra. The initial functional
can be regarded as a vacuum expectation value (functional) or a {\em state}
on $A$ if $A$ has a unit element.

A natural question which has been brought to our attention by 
K. Fredenhagen
is the following: can the GNS construction be extended to the associative
algebras occurring in deformation quantization? 
At first sight, this seems
almost impossible: firstly, for any {\em complex-valued} linear functional 
on $C^\infty(M)[[\lambda]]$
one would at once have an obvious convergence problem. Secondly,
although the deformed algebra does have an antilinear involution it does
not seem to share the remaining defining analytic properties of a complex
$C^*$-algebra such as the existence of a $C^*$-norm or any obvious uniform
structure in which it is complete.
This picture, however, will change drastically --as we should like to show
in the sequel-- when the linear functionals take their values not only 
in the
field of complex numbers but can have values in (some suitable field 
extension of) the ring $\mathbb C[[\lambda]]$ 
and are not only linear with respect to 
$\mathbb C$, but $\mathbb C[[\lambda]]$-linear:
This makes sense since
$C^\infty(M)[[\lambda]]$ is a $\mathbb C[[\lambda]]$-module 
in a natural way , a fact which lies at the heart of algebraic deformation
theory (see e.g. \cite {GS88,DL88}).
Thereby one stays in the formal category which avoids the
convergence problems. Moreover, in the subring $\mathbb R[[\lambda]]$
of formal power series with real coefficients there is an algebraic sense
of `asymptotic positivity': a power series in $\mathbb R[[\lambda]]$ is
defined to be positive (resp. negative) if and only if among its 
nonvanishing 
coefficients the one with lowest order is a positive (resp. negative)
real number. This defines a ring {\em ordering} (which is preserved 
under sums and
products). This structure allows us to speak of a positive linear 
functional
$\omega$ on the deformed algebra if and only if $\omega(\cc{f}*f)$ is 
nonnegative in $\mathbb R[[\lambda]]$ for all $f\in C^\infty(M)[[\lambda]]$.
For technical reasons it will be advantageous to replace the ring 
$\mathbb R[[\lambda]]$ by its quotient field $\LR$ of formal
Laurent series (with finite principal part) and likewise for
$\mathbb C[[\lambda]]$ and $C^\infty(M)[[\lambda]]$.
The traces considered in deformation quantization 
(see \cite [p. 3]{CFS92}; \cite [p. 229] {NT95a}, 
\cite [p. 155]{NT95b} where the Laurent field is denoted by 
$\mathbb C[\hbar^{-1}, \hbar]])$; and \cite [p. 171] {Fed96}) 
are particular examples for $\LC$-linear functionals.

By this simple replacement of the complex numbers by (field extensions of)
formal power series we get the following principal results:
\begin{itemize}
\item With the above notion of positive linear functionals with values
      in (a field extension of) $\LC$ we can imitate the
      algebraic part of the classical GNS construction for the algebras in
      deformation quantization and arrive at a representation of these
      algebras in a pre-Hilbert space over (a suitable field extension) of
      $\LC$ (section \ref {OrderSec}).
\item The above `asymptotic ordering' of the field defines a metric 
      topology on the field which can be used to define a norm 
      on the above pre-Hilbert space with values in 
      (a suitable extension of) $\LR$. This norm serves to define 
      a Cauchy completion of the pre-Hilbert space 
      (appendix \ref {AppA}, \ref{AppB}).
\item Apart from this algebraic construction this `formal' GNS 
      construction gives rise to the geometric interpretation 
      of these positive linear functionals as
      a {\em deformation (quantization) of classical positive linear 
      functionals}: these classical `states' in classical mechanics 
      are often given by certain measures (e.g. the Boltzmann 
      distribution) having support on certain submanifolds of $M$. 
      In this paper we shall deal with the simplest examples: 
      a) Dirac measures (whose support is one point in $M$) 
      for which we shall show in the case of the standard
      Wick star product in $\mathbb C^n$ that the above formal
      GNS construction gives rise to the formal
      Wick quantization with the correct formal Bargmann Hilbert space
      (section \ref {WickSec}), and
      b) the integral over configuration space $\mathbb R^n$ in 
      $\mathbb R^{2n}$ for the standard star product of Weyl type which
      will give rise to the usual symmetrization rule and the correct 
      formal $L^2$-space of (certain) square-integrable functions 
      of Schr\"odinger quantization on configuration space 
      (section \ref{ExampleIISec}). 
\item Moreover, we show that for the Fedosov star
      product of Wick type on an arbitrary K\"ahler manifold 
      (see \cite{BW96a} for details) the Dirac delta 
      functionals are positive functionals on the deformed 
      algebra which allows a representation of
      Bargmann type. The geometry of the K\"ahler manifold enters this
      representation by mapping the function $f\in \LS{C^\infty(M)}$
      on its Fedosov-Taylor series at the support of the delta functional
      (section \ref {KaehlerSec}).
\item Finally, in all the above examples we could also find a solution
      to the problem of convergence of formal powers series: either by
      using a given formal Hilbert base (in the K\"ahler examples) or a 
      set of suitable formal linear functionals (in the case of the 
      formal Schr\"odinger representation) we define a complex 
      Hilbert space out of the formal GNS Hilbert space for the 
      value $\lambda=\hbar$ in the following way: we first single 
      out the subset $H(\hbar)$ of all those elements $\psi$ of
      the formal Hilbert space which satisfy an infinite number 
      of suitable convergence conditions for all those formal 
      series with complex coefficients which arise when all the 
      scalar products of $\psi$ with the elements of the Hilbert 
      base or the values of all the given linear functionals on $\psi$ 
      are considered. This turns out to be a complex vector space and 
      the quotient space ${\cal H}(\hbar):=H(\hbar)/N(\hbar)$ 
      becomes a complex Hilbert space where $N(\hbar)$
      is the subspace of those elements of $H(\hbar)$ for which 
      all the above series vanish when $\lambda$ is replaced by 
      the real number $\hbar$. A certain complex subspace of the 
      formal algebra will respect this construction and can thus 
      be represented on (a subspace of) ${\cal H}(\hbar)$.

      These examples seem to support the point of view that the 
      convergence problem in deformation quantization should 
      be treated only {\em after} a formal GNS representation has 
      been chosen (by the choice of a suitable positive linear 
      functional), i.e. one should ``stay in the formal 
      category as long as possible''.
    
\end {itemize}
All these exemplary results seem to suggest that the 
concept of formal GNS construction sketched above may be a suitable 
candidate to formulate reasonable `Hilbert space' 
representations of the observable algebra within the framework 
of deformation quantization in a uniform manner:
the observable algebra is the primary object 
and the representing GNS (pre-) Hilbert spaces are subordinate 
objects parametrized by (often geometrically motivated) 
positive linear functionals (as opposed to other approaches to 
quantization where the observable algebra is defined as an 
operator algebra on a given Hilbert space).

The paper is organized as follows: Section \ref {MotivSec} 
recalls the necessary
notation and definitions in deformation quantization, introduces the use
of the field of formal Laurent series and gives some first elementary
examples and counterexamples of formal positive functionals. The main
flaw of the field of formal Laurent series is the fact that it is not
algebraically closed. 
Section \ref {OrderSec} is a bit technical: 
here we first describe the properties of general 
ordered fields and define pre-Hilbert spaces over such fields. 
This enables us to define the GNS pre-construction for associative 
algebras with involution on pre-Hilbert spaces. Further properties as 
topologies and absolute values of ordered fields including their 
algebraic closures and Cauchy completions are discussed in appendix 
\ref {AppA}. In appendix \ref {AppB} we finally describe Hilbert 
spaces and the general GNS construction.
The ideas of the proofs of many
of these results can be found in algebra text books, or in the literature
on $p$-adic functional analysis (see e.g. \cite{NBB71}), or in the works
on (nonclassical) Hilbert spaces (see e.g. \cite{GK77}, \cite{Kell80}, 
\cite{Gro79}).
In section \ref {GNSDeformSec} we apply the
results of section \ref {OrderSec} and the appendices 
\ref {AppA} and \ref{AppB} to the concrete fields we need 
for the GNS construction in deformation quantization, i.e. 
the algebraic closure of the Laurent
field, the so-called field of Newton-Puiseux series, and its Cauchy 
completion. The GNS construction is applied to the formal algebra
$\CNP{C^\infty(M)}$ of completed Newton-Puiseux series with coefficients
in $C^\infty(M)$ which contains the original deformed algebra
$C^\infty(M)[[\lambda]]$ as a subring. 
In section \ref {TimeSec} we prove that a formal
solution to the time-dependent Schr\"odinger equation exists 
as a well-defined curve of vectors in the formal GNS Hilbert 
space in case the classical Hamiltonian does not contain any 
negative powers of $\lambda$ and is contained in the Gel'fand 
ideal of the considered positive linear functional $\omega$.
Section \ref {WickSec} contains the above-mentioned
example of a GNS construction by means of the delta functional with
support at the origin for
the Wick star product in $\mathbb C^n$. In addition we obtain the
correct spectrum of the harmonic oscillator by an entirely formal
deduction. Section \ref {KaehlerSec} generalizes this
to the delta functional of an arbitrary point of an arbitrary K\"ahler
manifold equiped with the Fedosov star product of Wick type. 
In section \ref {ExampleIISec}
the above-mentioned example of the formal Schr\"odinger representation
in $\mathbb R^n$ is dealt with. Section \ref {OutSec} 
is an outlook on some open problems related to this approach.

\section {Motivation and basic concepts}
\label {MotivSec}

In this section we shall give a first heuristic motivation how one
can define positive linear functionals for a star product algebra 
without leaving the formal description.
First we have to introduce some notation:
Let $M$ be a $2n$-dimensional symplectic manifold.
The observable algebra in deformation quantization is given by
the formal power series $C^\infty (M)[[\lambda]]$
in the formal parameter $\lambda$ with coefficients in the smooth
complex-valued functions $C^\infty (M)$.
Then the addition in
$C^\infty (M)[[\lambda]]$ is the pointwise addition order by order in
$\lambda$ of the functions and the multiplication is a star product for
$M$. A star product of two functions $f, g \in C^\infty (M)$ is a formal
power series
\[
    f * g = \sum_{r=0}^\infty \lambda^r M_r (f, g)
\]
with bilinear operators $M_r:C^\infty(M)\times C^\infty(M)\rightarrow
C^\infty(M)$ such that $*$ extends to a
$\mathbb C[[\lambda]]$-bilinear associative product for
$C^\infty (M)[[\lambda]]$ with the following properties:
The lowest order is the pointwise multiplication,
i.e. $M_0 (f,g) = fg$, and the first order commutator is given by
$i\lambda$ times the Poisson bracket, i. e.
$M_1(f,g) - M_1 (g,f) = i\{f,g\}$.
The constant function $\bfmath 1$ is the unit element with respect to $*$,
i.e. $f * \bfmath 1 = \bfmath 1 * f = f$.
We assume for simplicity that all the operators $M_r$ are not only local but
bidifferential operators. We demand in addition that the complex
conjugation of elements in $C^\infty (M)[[\lambda]]$
(where $\cc \lambda := \lambda$) is an antilinear algebra involution:
\[ 
    \cc {f*g} = \cc g * \cc f
\]
Then $C^\infty (M)[[\lambda]]$ together with a star product is not only
a $^*$-algebra over the field $\mathbb C$ but also a
$\mathbb C[[\lambda]]$-module. Note that $\lambda$ is to be 
identified directly with $\hbar$.

It is now an easy step to pass form the
ring of formal power series to the corresponding quotient field which
turns out to be the field of {\em formal Laurent series}. 
The field of formal Laurent series with real (and analogously
complex) coefficients is defined by
\BEQ {RealLaurent}
    \LR := \left.\left\{ \sum_{r=-N}^\infty \lambda^r a_r
    \; \right| \; a_r \in \mathbb R, N \in \mathbb Z \right\} .
\EEQ
Then $\LR$ has the structure of a field if we define the addition order
by order in $\lambda$ and the multiplication by
\BEQ {LMultDef}
     \left(\sum_{r=-N}^\infty \lambda^r a_r\right)
     \left(\sum_{s=-M}^\infty \lambda^s b_s\right)
     :=
     \sum_{t=-N-M}^\infty \lambda^t
     \sum_{r+s=t \atop r \ge -N, s \ge -M} a_rb_s .
\EEQ
Then one easily shows that $\LR$ and $\LC$ are fields and with
$\cc\lambda := \lambda$ we have the following natural inclusions:
$\mathbb R \subset \mathbb R[[\lambda]] \subset \LR$,
$\mathbb C \subset \mathbb C[[\lambda]] \subset \LC$ and
$\LR \subset \LC$. Moreover $\LR$ is the quotient field of the ring
$\mathbb R[[\lambda]]$ and now $\lambda$ has an inverse,
namely $\lambda^{-1}$.

Remark: From the physical point of view it will 
be necessary to consider negative powers of Planck's constant anyway
since $\hbar$ will appear in the denominator in important situations
as for example in the energy eigenvalues of the hydrogen atom.

More generally, if $V$ is a $\mathbb K$-vector space
($\mathbb K=\mathbb R,\mathbb C$) then we define the
corresponding vector space of formal Laurent series $\LS V$ by
\BEQ {FormalLaurent}
    \LS V := \left.\left\{ \sum_{r=-N}^\infty \lambda^r v_r
    \; \right| \; v_r \in V, N \in \mathbb Z \right\}
\EEQ
and notice that $\LS V$ is a $\LS{\mathbb K}$-vector space in the obvious way.
Moreover the field $\LR$ is known to be an ordered field with a unique
ordering (e. g. \cite [p. 73]{Rui93}):
\begin {lemma} \label {OrderLem}
$\LR$ is an ordered field with a unique non-archimedian ordering
relation such that $\lambda > 0$. The positive elements are given by
\BEQ {LRPosDef}
    \sum_{r=-N}^\infty \lambda^r a_r > 0
    \iff
    a_{-N} > 0 .
\EEQ
\end {lemma}
A symbolic picture of this non-archimedian ordering is given by
\[
   \cdots < \lambda^{-1} \mathbb  R^-
   < \mathbb  R^-
   < \lambda \mathbb  R^-
   < \cdots
   < 0
   < \cdots
   < \lambda \mathbb  R^+
   < \mathbb  R^+
   < \lambda^{-1} \mathbb  R^+
   < \cdots
   .
\]
Apart from this ordering relation
the fields $\LR$ and $\LC$ allow an {\em absolute value}
$\varphi: \LC \to \mathbb R$ defined for $a \in \LC$ by
\[
    \varphi (a) := 2^{-o(a)}
\]
where $o(a) \in \mathbb Z$ is the order of the
lowest non-vanishing term in the formal Laurent series $a$.
Both structures, the order and the absolute value, 
can be used to define {\em topologies} on the field $\LR$
(see appendix \ref {AppA} for definitions)
which will turn out to be the same 
(see section \ref {GNSDeformSec}). An important
observation will be that in $\LR$ every Cauchy 
sequence defined with respect to the order converges 
in $\LR$.

Using the fields $\LR$ and $\LC$ instead of $\mathbb R$ and $\mathbb C$
we can generalize the $\mathbb C$-algebra of observables
$C^\infty (M)[[\lambda]]$ in the obvious way to a $\LC$-algebra
$\LS{C^\infty (M)}$ where we extend the star product to a $\LC$-bilinear
product of $\LS{C^\infty (M)}$. Now we can consider $\LC$-linear
functionals $\omega: \LS{C^\infty (M)} \to \LC$ and define positive
linear functionals using the complex conjugation as antilinear
algebra involution and the ordering of $\LR$ completely analogously
to the case of $C^*$-algebras over $\mathbb C$ \cite {BR87}:
\begin {definition} \label {PosStateDef}
A $\LC$-linear functional $\omega: \LS{C^\infty (M)} \to \LC$ is
called positive iff
\[
    \omega (\cc f * f) \ge 0 \qquad \forall f \in \LS{C^\infty (M)}
\]
and a state iff $\omega$ is positive and $\omega (\bfmath 1) = 1$.
\end {definition}
We will now give simple examples for positive linear functionals.
Using the locality of the star product we can prove the following
lemma:
\begin {lemma}
Let $0 \le \varrho \in C^0 (M)$ with compact support. Then
$\omega_\varrho$ defined by
\[
    \omega_\varrho (f) := \int_M f\varrho \, \Omega
    \qquad \qquad
    (\Omega \mbox{ symplectic volume form})
\]
is a positive linear functional of $\LS{C^\infty (M)}$ with respect to
{\em any} star product and
$\omega_\varrho (\cc f * f) = 0 \iff f\varrho \equiv 0$.
\end {lemma}
A counterexample is given by the {\em delta-functionals}
in $\mathbb R^2$ and the Weyl-Moyal star product (\ref{WeylProdDef}). 
We easily compute
that the evaluation functional $\delta_{(q_0,p_0)}$ at
$(q_0, p_0) \in \mathbb R^2$ is not positive since
\[
    \delta_{(q_0, p_0)} \left(
    \cc{\left( (q-q_0)^2 + (p-p_0)^2\right)} *
    \left( (q-q_0)^2 + (p-p_0)^2\right)\right)
    = -\frac{1}{2} \lambda^2 < 0
\]
in spite of the fact that the delta functional is a positive
functional with respect to the pointwise product and hence a classical
state.

Now we can construct a representation of the
$\LC$-algebra $\LS{C^\infty (M)}$ using a given positive linear
functional in the same way as the usual GNS representation of a
$C^*$-algebra is constructed. But this will be done in a
more general context in the next section and then it will turn out 
that there is even a better choice than the formal Laurent series 
since the field $\LR$ is not {\em real closed} and $\LC$ is not 
{\em algebraically closed} since for example the positive element 
$\lambda$ is no square.

\section {Ordered fields and general GNS pre-construction}
\label {OrderSec}

In this rather technical section we shall present some general
properties of ordered fields and
vector spaces and algebras over such fields.
All statements can in principal be proved in a fashion very similar
to the case of real or complex numbers,
and for the proofs only the ordering axioms will be used.
In order to make this exposition not too long we shall omit most of the
proofs. As general references we have used the textbooks by Jacobson
\cite{JacI} (in particular Chapter 5) and \cite{JacII} (Chapters 9, 11),
by Ruiz \cite{Rui93}, Kelley \cite{Kel55}, Rudin \cite{Rud87}, Yosida
\cite{Yos80}, Bratteli-Robinson \cite{BR87}, and Haag \cite{Haa93}.

\begin {definition} [{\cite [Def. 5.1] {JacI}}] \label {OrderDef}
An ordered field $(\field R, P)$ is field $\field R$ with a subset
$P$ of {\em positive elements} of $\field R$ such that $0 \not\in P$ and if
$a \in \field R$ then either $a=0$, $a \in P$ or $-a \in P$ and if
$a,b \in P$ then $a+b \in P$ and $ab \in P$. Then an ordering relation
is defined by $a<b$ iff $b-a \in P$.
\end {definition}
The symbols $>$, $\le$ and $\ge$ will be used as in the case
$\field R = \mathbb R$.
The elements with $a<0$ will be called negative.
The set of positive elements will also be denoted by $\field R^+$ and
the negative elements by $\field R^-$.
Clearly, every square $a^2$ is positive or $0$ and hence $-1<0<1$.
Since $n1 = 1 + \cdots + 1$ is a sum of positive numbers for all 
$n \in \mathbb N$ we have $n1 > 0$ and hence $\field R$
has characteristic zero.
We define $|a| := a$ if $a \ge 0$ and $|a| := -a$ if
$a<0$. Then $|\cdot|$ clearly satisfies
\BEQ {AbsProps}
     \begin {array} {c}
     |a| \ge 0 \quad {\rm and} \quad |a| = 0 \Longleftrightarrow a = 0 \\
     |ab| = |a||b| \\
     \big| |a| - |b| \big| \le |a+b| \le |a| + |b| .
     \end {array}
\EEQ  
The ordering relation $<$ is called {\em archimedian} iff there is for any 
$0<a$ and $b \in \field R$ a natural number $n \in \mathbb N$ such that 
$na > b$. Otherwise the ordering will be called {\em non-archimedian}.

Now we fix once and for all the quadratic field extension
$\field C := \field R(i) (=\field R\oplus i \field R)$ of an ordered field
$\field R$ where $i^2 := -1$.
In $\field C$ we define complex conjugation as usual by
$a + ib \in \field C \; \mapsto \; \cc{a+ib} := a - ib$
where $a,b \in \field R$.
Then the complex conjugation is an involutive field automorphism.
The elements of the subfield $\field R$ of $\field C$ will be called
{\em real} and are characterised as usual by
$a \in \field R \subset \field C \Longleftrightarrow \cc a = a$.
Clearly $\cc a a \in \field R$, $\cc a a \ge 0$
and $\cc a a = 0 \Longleftrightarrow a = 0$
for any $a \in \field C$ and we define $|a|^2 := \cc a a$.

Now we transfer definitions and some simple algebraic
results from the theory of complex vector spaces with Hermitian products
to the general case of a vector space over an arbitrary ordered
field $\field R$ and its quadratic field extension $\field C$.
\begin {definition} \label {preHilbertDef}
Let $\field C = \field R(i)$ be the quadratic
field extension of an ordered field $\field R$  and 
$\HC$ a $\field C$-vector space.
A map $\SP\cdot\cdot :
\HC \times \HC \to \field C$
is called a Hermitian product iff for all
$\phi,\psi,\chi \in \HC$ and for all
$a,b \in \field C$
\begin {enumerate}
\item $\SP \cdot\cdot$ is antilinear in the first argument, i.e.
      $\SP{a\phi + b\psi}{\chi} =
      \cc a \SP{\phi}{\chi} + \cc b \SP{\psi}{\chi}$
\item $\SP \phi\psi= \cc {\SP \psi\phi}$ 
\item $\SP \phi\phi \ge 0$ and
      $\SP \phi\phi = 0 \; \Longleftrightarrow \; \phi=0$ .
\end {enumerate}
A $\field C$-vector space with a Hermitian product $\SP \cdot\cdot$
is called a pre-Hilbert space. A linear map
$U: \HC \to \KC$ from one
pre-Hilbert space over $\field C$ to another is called an isometry iff
$\SP {U\phi}{U\psi} = \SP \phi\psi$
for all $\phi,\psi \in \HC$ and unitary iff $U$ is a bijective isometry.
\end {definition}

For a pre-Hilbert space we have the
Cauchy-Schwarz inequality for the Hermitian product which is a simple
imitation of the standard arguments over the complex numbers
(cf e.g. \cite[p. 77]{Rud87}):
\begin {lemma} [Cauchy-Schwarz inequality]
Let $\HC$ be a pre-Hilbert space with Hermitian
product $\SP \cdot\cdot$ and $\phi,\psi \in \HC$. Then
\BEQ {CSI}
    \SP \phi\psi \SP \psi\phi \; \le \; \SP \phi\phi \SP \psi\psi
\EEQ
with equality if and only if $\psi$ and $\phi$ are linearly dependent.
\end {lemma}

Furthermore, for a linear operator from
$\HC$ into $\HC$ we want to define an adjoint operator.
In general it is rather complicated to see whether an adjoint operator
exists or not. Without further assumptions we can only state the
following definition for everywhere defined linear operators:
\begin {definition} [Adjoint operator] \label {AdjointDef}
Let $\HC$ be a pre-Hilbert space and
$A : \HC \to \HC$ a linear map. Then a linear map
$B : \HC \to \HC$ is called an adjoint operator to $A$ iff
$\SP {\phi}{A\psi} = \SP{B\phi}{\psi}$ for all $\phi, \psi \in \HC$.
In this case we write $B = A^*$.
\end {definition}
\begin {lemma}
Let $\HC$ be a pre-Hilbert space and
$A,B: \HC \to \HC$ linear operators.
If $A^*$ exists then it is unique and if $A^*$ and $B^*$ exist
then $(aA+bB)^*$, $(A^*)^*$ and $(AB)^*$ exist for
$a,b \in \field C$ and $(aA + bB)^* = \cc a A^* + \cc b B^*$,
$(A^*)^* = A$ and $(AB)^* = B^*A^*$.
\end {lemma}

Let $\mathcal A$ be an associative and not necessarily commutative
algebra over the quadratic field extension
$\field C = \field R(i)$ of
an ordered field $\field R$.
We shall only consider algebras with an antilinear involution $^*$
compatible with the complex conjugation in $\field C$. 
This means that there is a map
$^* : \mathcal A \to \mathcal A$ such that for all $A,B \in \mathcal A$ 
and $a\in \field C$
\BEQ {AlgInvDef}
    \begin {array} {ccc}
    (A+B)^* = A^* + B^* & \quad & (aA)^* = \cc a A^* \\
    (AB)^* = B^* A^*    & \quad & (A^*)^* = A .
    \end {array}
\EEQ
If the algebra $\mathcal A$ has a unit element $\bfmath 1$ 
then necessarily
$\bfmath 1^* = \bfmath 1$. By means of the algebra involution $^*$
and the ordering relation in $\field R$ we can define
{\em positive linear functionals} analogously to the case of
$C^*$-algebras over the complex numbers $\mathbb C$ 
\cite {BR87}, \cite{Haa93}:
\begin {definition}
Let $\mathcal A$ be a $\field C$-algebra with involution $^*$ and 
$\omega : \mathcal A \to \field C$ a linear functional where
$\field C = \field R(i)$ and $\field R$ is an ordered
field. Then $\omega$
is called positive iff for all $A \in \mathcal A$
\[
    \omega (A^* A) \ge 0 .
\]
If in addition $\mathcal A$ has a unit element $\bfmath 1$ then
$\omega$ is called
a state iff $\omega$ is positive and $\omega (\bfmath 1) = 1$.
\end {definition}
For a positive linear functional we can prove the
Cauchy-Schwarz inequality:
\begin {lemma} [Cauchy-Schwarz inequality] \label {CSILem}
Let $\mathcal A$ be a $\field C$-algebra with involution $^*$ and 
$\omega : \mathcal A \to \field C$ a positive linear functional. 
Then we have for all $A, B \in \mathcal A$:
\begin {eqnarray}
    \omega (A^* B) & = & \cc {\omega (B^* A)}  \label {CSi} \\
    \omega (A^* B) \cc{\omega (A^* B)}
    & \le & \omega (A^* A) \omega (B^* B) \label {CSii}
\end {eqnarray}
Moreover $\omega (\bfmath 1) = 0$ implies $\omega \equiv 0$ and
$\omega (A^*) = \cc {\omega (A)}$ and if
$\omega (A^*A) = 0$ then $\omega (A^*B) = 0$ for all
$B \in \mathcal A$.
\end {lemma}

For an algebra $\mathcal A$ with involution $^*$ 
over the quadratic field extension $\field C$ of an 
ordered field $\field R$ we are able to transfer the well-known 
GNS construction of representations for $C^*$-algebras:
We construct a {\em representation} of $\mathcal A$ in a pre-Hilbert
space over $\field C$.
Firstly, we shall only deal with the algebraic properties of this
construction and examine analytic properties in appendix \ref {AppB}.
We are mainly using the notation of \cite {BR87}.
Let $\omega : \mathcal A \to \field C$ be a positive linear functional.
Then we consider the following subspace of $\mathcal A$:
\BEQ {Gelfandideal}
    \mathcal J_\omega \; := \; \{ A \in \mathcal A \; | \;
    \omega (A^*A) = 0\}
\EEQ
This subspace $\mathcal J_\omega$ is called the {\em Gel'fand ideal}
of $\omega$ and by means of lemma \ref {CSILem} it is easily proved that
$\mathcal J_\omega$ is indeed a 
left ideal of $\mathcal A$. In the next step one considers the 
quotient vector space
\BEQ {AquotientJ}
    \mathfrak H_\omega \; := \; \mathcal A / \mathcal J_\omega 
\EEQ
where the equivalence classes in $\mathfrak H_\omega$ are denoted by
\BEQ {Hvectors}
    \psi_A \; := \; \{A' \in \mathcal A \; | \; A' = A + I,
               I \in \mathcal J_\omega\} .
\EEQ
On the quotient space $\mathfrak H_\omega$ we can define a Hermitian
product by
\BEQ {hermprod}
    \SP {\psi_A}{\psi_B} \; := \; \omega (A^*B)
\EEQ
which is well-defined since $\mathcal J_\omega$ is a left ideal. 
Furthermore this definition leads indeed to a {\em non-degenerate} Hermitian
product for $\mathfrak H_\omega$ which will make $\mathfrak H_\omega$
a pre-Hilbert space over $\field C$.
In a third step one defines a representation $\pi_\omega$
of $\mathcal A$ on $\mathfrak H_\omega$ by
\BEQ {piOmegaRep}
    \pi_\omega (A) \psi_B \; := \; \psi_{AB} .
\EEQ
To prove that this is well-defined we need once more the fact that
$\mathcal J_\omega$ is a left ideal. 
Then we notice that $\pi_\omega$ is a $^*$-representation: 
$\pi_\omega (AB) = \pi_\omega (A)\pi_\omega (B)$, $\pi_\omega$ is linear
and $\pi_\omega (A^*) = \pi_\omega (A)^*$.
In this case it is checked directly that $\pi_\omega (A^*)$ is 
the adjoint operator of $\pi_\omega (A)$ in the sense of 
definition \ref{AdjointDef}.
If $\mathcal A$ has a unit element $\bfmath 1$ and $\omega$ is a
state then the 
representation $\pi_\omega$ is {\em cyclic} with the cyclic (vacuum) vector
$\psi_{\bfmath 1}$ since every vector $\psi_A \in \HC$
can be written as $\psi_A = \pi_\omega (A) \psi_{\bfmath 1}$ and we have
\BEQ {CyclOmega}
    \omega (A) = \SP{\psi_{\bfmath 1}}{\pi_\omega (A) \psi_{\bfmath 1}}.
\EEQ
In this case we show analogously to the case of $C^*$-algebras
that the representation $\pi_\omega$ is unique up
to unitary equivalence:
If $\mathfrak H', \pi', \psi'$ is another cyclic 
$^*$-representation with cyclic unit vector $\psi'$ such that
$\omega (A) \; = \; \SP{\psi'}{\pi'(A)\psi'}$
for all $A \in \mathcal A$ then we have a unitary map 
$U: \mathfrak H_\omega \to \mathfrak H'$ such that
$U^{-1} \pi'(A) U  = \pi (A)$ and $U\psi_{\bfmath 1} = \psi'$.
We shall resume this in the following proposition:
\begin {proposition} [General GNS pre-construction]
\label {GNSConst}
Let $\mathcal A$ be a $\field C$-algebra with involution $^*$
where $\field C = \field R(i)$
is the quadratic field extension of an ordered field $\field R$.
For any positive linear functional there exists a $^*$-representation
$\pi_\omega$ on a pre-Hilbert space $\mathfrak H_\omega$
as constructed above which is called the 
{\em GNS representation on $\mathfrak H_\omega$}.
If $\mathcal A$ in addition has a unit element $\bfmath 1$ and
$\omega$ is a state then this representation is cyclic and we have
$\omega (A) = \SP{\psi_{\bfmath 1}}{\pi_\omega (A) \psi_{\bfmath 1}}$
and this property defines this representation up to unitary equivalence.
\end {proposition}
The following obvious generalization will be very useful if the positive
linear functional is only defined on a proper ideal of $\mathcal A$,
a situation which will occur in section 7.
\begin {corollary}
Let $\mathcal A$ be a $\field C$-algebra with involution $^*$ and
$\mathcal B \subset \mathcal A$ a twosided ideal of $\mathcal A$ with
$\mathcal B = \mathcal B^*$. Let $\omega : \mathcal B \to \field C$ 
be a positive linear functional and denote by $\pi_\omega$ the
GNS representation of $\mathcal B$ on 
$\mathfrak H_\omega = \mathcal B / \mathcal J_\omega$ where 
$\mathcal J_\omega \subset \mathcal B$ is the Gel'fand ideal 
of $\omega$. Then $\mathcal J_\omega$ is a left ideal of $\mathcal A$ 
and $\pi_\omega$ can be extended canonically to a representation 
of $\mathcal A$ on $\mathfrak H_\omega$.
\end {corollary}
\begin {proof}
Let $a \in \mathcal A$ and $b \in \mathcal J_\omega$. Since $B$ is an ideal
it follows that
$ab$ and $a^*ab$ are contained in $B$. Moreover 
$\omega((ab)^*ab)\cc{\omega((ab)^*ab)} =
 \omega(b^*(a^*ab))\cc{\omega(b^*(a^*ab))}$
$\leq  \omega(b^*b)\omega((a^*ab)^*a^*ab)=0$ 
thanks to the Cauchy Schwarz
inequality which implies that $ab\in J_\omega$. The rest of the Corollary
follows easily.    
\end {proof}

\section {The GNS Construction in deformation quantization}
\label {GNSDeformSec}

Since the Laurent field $\LR$ (definition (\ref {RealLaurent}))
is not real closed we are searching for a field extension of $\LR$
such that this extension will be real closed {\em and}
Cauchy complete with respect to the ordering relation since
this is the important property we need if we want to consider Hilbert
spaces as in appendix \ref {AppB}.

First we define the following systems of subsets of the 
rational numbers $\mathbb Q$:
\begin {definition} \label {SuppSysDef}
$S \subset \mathbb Q$ is called NP-admissible iff $S$ has a
smallest element and there is a positive integer such that
$N\cdot S \subset \mathbb Z$ and $S$ is called CNP-admissible
iff $S$ has a smallest element and $S \cap [i,j]$ is finite
for any $i,j \in \mathbb Q$. Then we define
$\mathcal S_{\rm CNP} := \{ S \in \mathbb Q \; | \; 
S \mbox{ \rm is CNP-admissible} \}$ and analogously
$\mathcal S_{\rm NP}$.
\end {definition}
Let $V$ be a vector space over $\mathbb R$ (or $\mathbb C$) and
$f: \mathbb Q \to V$. Then we define the support of $f$ by
$\supp f := \left\{ q \in \mathbb Q \; \big| \; f(q) \ne 0 \right\}$.
\begin {definition}
Let $V$ be a vector space over $\mathbb R$ (or $\mathbb C$).
Then we define the formal Newton-Puiseux series (NP series) with 
coefficients in $V$ by
\BEQ {NPVR}
    \NP V := \left\{ f: \mathbb Q \to V \; \big |\;
             \supp f \in \mathcal S_{\rm NP} \right\}
\EEQ
and the completed Newton-Puiseux series (CNP series) with coefficients
in $V$ by
\BEQ {CNPVR}
    \CNP V := \left\{ f: \mathbb Q \to V \; \big |\;
              \supp f \in \mathcal S_{\rm CNP} \right\} .
\EEQ
\end {definition}
Clearly $\NP V$ and $\CNP V$ are both vector spaces 
over $\mathbb R$ (or $\mathbb C$) and elements in $\NP V$ and $\CNP V$ are
written in the form
\[
    f = \sum_{q\in \supp f} \lambda^q f_q \qquad f_q := f(q) .
\]
and since $\mathcal S_{\rm NP} \subset \mathcal S_{\rm CNP}$ we notice
that $V \subset \LS V \subset \NP V \subset \CNP V$.
If $f \in \NP V$ then either the support is finite or it can be written
as $f = \sum_{r=-M}^\infty \lambda^{\frac{r}{N}} f_r$ where
$M \in \mathbb Z$ and $N \in \mathbb N$. If $f \in \CNP V$ then the
support is again either finite or can be written as a sequence
$q_0 < q_1 < \cdots \in \mathbb Q$ without any accumulation points.
In both cases we will write
\[
    f = \sum_{r=0}^\infty \lambda^{q_r} f_{q_r}.
\]
For $f \in \CNP V$
we define the order $o(f) := \min(\supp f)$ for $f \ne 0$ and set
$o(0) := +\infty$ and we define
\BEQ {CNPAbsValDef}
    \varphi (f) := 2^{-o(f)} \mbox { for } f \ne 0
    \qquad  \varphi(0) := 0
\EEQ
which leads to an {\em ultra-metric}
$d_\varphi (f, g) := \varphi (f-g)$ for the vector space $\CNP V$.
\begin {proposition} \label {CNPCompProp}
Let $V$ be a vector space over $\mathbb R$ (or $\mathbb C$). Then
$\CNP V$ is a complete metric space with respect to the
metric $d_\varphi$ and $\NP V$ is dense in $\CNP V$.
\end {proposition}
\begin {proof}
Let $(f^{(n)})_{n \in \mathbb N}$ be a Cauchy sequence
in $\CNP V$ with respect to the metric $d_\varphi$. Then
there are natural numbers $N_0 \le N_1 \le N_2 \le \cdots \in \mathbb N$
such that $d_\varphi (f^{(n)}, f^{(m)}) < 2^{-k}$ for all $n,m \ge N_k$.
Then we define $f_q := f^{(N_0)}_q$ for $-\infty < q \le 0$ and
$f_q := f^{(N_k)}_q$ for $k-1 < q \le k$ where $k \in \mathbb N$.
Then $f = \sum_q \lambda^q f_q$ is a well-defined series in
$\CNP V$ since $\supp f$ has a smallest element and in any
intervall of the form $[i,i+1]$ with $i \in \mathbb Z$ are only
finitely many $f_q \ne 0$ for $q \in [i,i+1]$ and hence
$\supp f \in \mathcal S_{\rm CNP}$. Since the $f^{(n)}$ coincide
with $f$ up to a sufficiently increasing order thanks to the
Cauchy condition it follows that $f^{(n)} \to f$ which implies that
the metric space $\CNP V$ is complete.
Now let $f = \sum_q \lambda^q f_q$ be an arbitrary element in
$\CNP V$. Then every element of the sequence
$f^{(m)} := \sum_{q \le m} \lambda^q f_q$ where $m \in \mathbb N$
has finite support and hence $f^{(m)} \in \NP V$.
Moreover clearly $f^{(m)} \to f$ which proves the proposition.
\end {proof}

\begin {corollary}
Let $V$ be a vector space over $\mathbb R$ (or $\mathbb C$) and
$f = \sum_q \lambda^q f_q \in \CNP V$. Then
$f^{(m)} := \sum_{q \le m} \lambda^q f_q$ converges to $f$.
\end {corollary}

We will now concentrate on the vector space $\NPR$ and $\CNPR$ (and
analogously for $\NPC$ and $\CNPC$) itself.
Analogously to the case of formal Laurent series we define a
multiplication of two elements $a,b \in \CNPR$ by
\BEQ {CNPMultDef}
    ab = \left(\sum_{p \in \supp a} \lambda^p a_p \right)
         \left(\sum_{q \in \supp b} \lambda^q b_q\right)
    := \sum_{t \in \supp a + \supp b}
       \lambda^t \sum_{p+q = t} a_pb_q
\EEQ
where $\supp ab = \supp a + \supp b := \{ p+q \; | \; p \in \supp a, q
\in \supp b \}$. Then $\supp ab \in \mathcal S_{\rm CNP}$ and
in any order $t \in \supp ab$ the sum $p+q = t$ for $p \in \supp a$ and
$q \in \supp b$ is finite. Hence $ab$ is again a well-defined CNP
series. A proof in a more general context could be found
in \cite [p. 81] {NBB71}.
Moreover the product is clearly associative and commutative and we 
can find
for any $0 \ne a \in \CNPR$ an inverse $a^{-1} \in \CNPR$ and hence
$\CNPR$ becomes a field. Furthermore one can show that
$\NPR$ is a subfield of $\CNPR$. A similar result holds for $\CNPC$.

If $V$ is a vector space over $\mathbb K$
($\mathbb K = \mathbb R$ or $\mathbb C$)
then $\NP V$ is a vector space over $\NP{\mathbb K}$ and $\CNP V$ is a
vector space over $\CNP{\mathbb K}$. We call these vector spaces the
{\em canonical extensions} of $V$ to vector spaces over the NP series
and the CNP series with real (or complex) coefficients.
Let $V^*$ be the (algebraic) dual vector space of $V$ then any
$\mathbb K$-linear functional $\omega \in V^*$ has an obvious
canonical $\CNP{\mathbb K}$-linear extension to a
functional of $\CNP V$ by applying $\omega$ order by order.
This extension will always be understood.
A particular subspace of all $\CNP{\mathbb K}$-linear functionals
$\omega : \CNP V \to \CNPR$ is given by $\CNP{V^*}$,
i.e. the linear functionals of the form
$\omega = \sum_{q \in \supp \omega} \lambda^q \omega_q$ where
$\omega_q \in V^*$ and $\supp \omega \in \mathcal S_{\rm CNP}$.
Note that $\CNP{V^*}$ is in general a {\em proper} subspace of
all $\CNP{\mathbb K}$-linear functionals of $\CNP V$,
i.e. $\CNP{V^*} \subsetneq (\CNP V)^*$ iff $V$ is
infinite-dimensional.

The field of CNP series $\CNPR$ (or $\CNPC$)
has a canonical non-archimedian absolute value, namely $\varphi$ and
$\CNPR$ has a unique order which turns out to be compatible with the
absolute value:
\begin {proposition}
The map $\varphi$ defined as in
(\ref {CNPAbsValDef}) is a non-archimedian and non-trivial absolute
value for the fields $\CNPR$ resp. $\CNPC$.
The field $\CNPR$ is an ordered field and this unique order is
compatible with the absolute value $\varphi$ and hence
$\mathcal T = \mathcal T_\varphi$ and $\mathcal U = \mathcal U_\varphi$.
The positive elements are given by $a > 0$ iff $a_{q_0} > 0$ where
$q_0 = \min(\supp a)$.
\end {proposition}
\begin {theorem} \label {NPTheo}
\begin {enumerate}
\item $\NPC \cong \NPR (i) (=\NPR \oplus i \NPR)$ is algebraically
      closed and $\NPR$ is real closed.
\item $\CNPC \cong \CNPR (i) (= \CNPR \oplus i \CNPR)$ is 
      algebraically closed
      and $\CNPR$ is real closed and both are Cauchy complete with
      respect to the ordering relation of $\CNPR$.
\end {enumerate}
\end {theorem}
\begin {proof}
The first part is the Newton-Puiseux theorem (\cite [p. 61] {Rui93} or
\cite [p. 595] {JacII} for a more general case).
To prove the second part we notice that
$\CNPC$ is the completion of $\NPC$ with respect to the metric
$d_\varphi$ induced by the absolute value $\varphi$ according to
proposition \ref {CNPCompProp}. Hence
K\"ursch\'ak's theorem ensures that $\CNPC$ is again algebraically
closed, see \cite [p. 584] {JacII}.
Furthermore, since $\CNPC \cong \CNPR \oplus i \CNPR$ it follows e.g.
from the theorem of Artin and Schreier \cite [p. 674] {JacII} that
$\CNPR$ is real closed.
\end {proof}

This theorem ensures that the fields $\CNPR$ and $\CNPC$ have indeed
the algebraic and analytic properties needed for the GNS construction in
deformation quantization and the definition of Hilbert spaces
over $\CNPC$. Note that the NP series would not suffice since
they are not Cauchy complete. Using the field $\CNPC$
we consider the $\CNPC$-vector space $\CNP{C^\infty (M)}$ and extend the
star product to a $\CNPC$-bilinear product on $\CNP{C^\infty (M)}$.
Hence the observable algebra in deformation quantization becomes a
$\CNPC$-algebra. To define states and GNS representations we need
positive linear functionals of this algebra and first we notice that
the results of section \ref {MotivSec} not only hold for the field of
formal Laurent series but also for the CNP case.
Then we formulate the GNS construction as a corollary to theorem 
\ref {GNSConst} and \ref {FinalGNSConst}:
\begin {corollary} [GNS (pre-)construction in deformation quantization] 
{$\;$}
Let $M$ be a symplectic manifold and $*$ a star product for $M$ 
with $\cc{f*g} = \cc g * \cc f$. 
Let $\field K = \LC$ or $\NPC$ or $\CNPC$, and let 
$\mathcal A_{\field K} := \LS{C^\infty (M)}$, 
$\NP{C^\infty (M)}$, $\CNP{C^\infty (M)}$, respectively. 
Then for any positive linear 
functional $\omega: \mathcal A_{\field K} \to \field K$ 
there exists a pre-Hilbert space $\mathfrak H_\omega$ 
over $\field K$ 
carrying a $^*$-representation
$\pi_\omega$ of $\mathcal A_{\field K}$ 
as constructed in theorem \ref {GNSConst}.
We shall call this representation of $\mathcal A_{\field K}$ 
the {\em GNS representation on $\mathfrak H_\omega$}.

In the case $\field K = \CNPC$ the pre-Hilbert space can be
completed to a Hilbert space $\hat{\mathfrak H}_\omega$ 
over $\CNPC$ as in theorem \ref {FinalGNSConst}. Then we shall call 
$\pi_\omega$ the {\em GNS representation of $\CNP{C^\infty (M)}$ 
in $\hat{\mathfrak H}_\omega$}.
 
If $\omega$ is a state then the GNS representation is cyclic and we have 
$\omega (f) = \SP{\psi_{\bfmath 1}}{\pi_\omega (f) \psi_{\bfmath 1}}$
and this property defines this representation up to unitary equivalence.
\end {corollary}

If we consider positive functionals of the particular form
$\omega = \sum_{q \in \supp \omega} \lambda^q \omega_q$
where $\omega_q : \CNP{C^\infty (M)} \to \mathbb C$ and
$\supp \omega \in \mathcal S_{\rm CNP}$ as mentioned above then we
can show that for a {\em state} no negative powers in the
formal parameter could occur using the Cauchy-Schwarz inequality for
the state:
\begin {lemma}
Let $\omega: \CNP {C^\infty (M)} \to \CNPC$ be a positive
linear functional of the form
$\omega = \sum_{q \in \supp \omega} \lambda^q \omega_q$
such that $\omega(\bfmath 1) = 1$.
Then $q_0 := \min (\supp \omega) = 0$ and $\omega_0$ is
a classical state, i.e. $\omega_0 (\cc f f) \ge 0$ for
$f \in C^\infty (M)$.
\end {lemma}
\begin {proof}
Let $f \in \CNP{C^\infty (M)}$ be a function such that
$\omega_{q_0} (f) \ne 0$ where we can assume that $f \in C^\infty (M)$.
Then $o(\omega (f) \omega (\cc f)) = 2q_0$ and $o(\cc f *f)) \ge q_0$
which is in contradiction to the Cauchy-Schwarz inequality
$\omega (f) \omega (\cc f) \le \omega (\cc f * f)$ if $q_0 < 0$
where we used $\omega (\bfmath 1) = 1$. Moreover this implies
$\omega_{q_0} (\cc f f) \ge 0$ since $\cc f * f = \cc f f + \ldots$.
\end {proof}

This lemma has the following physical interpretation: Not only
the classical observable algebra $C^\infty (M)$ is deformed but also
{\em the classical states are deformed} to obtain the states for the
deformed algebra. Hence the restriction to these particular linear
functionals in $\CNP{(C^\infty(M)^*)}$ fits well into the general concept
of deformation. As we shall see in the examples the positive linear 
functionals are typically formal CNP series with coefficients even in a 
smaller space namely in $C^\infty (M)'$ resp. $C^\infty_0 (M)'$ 
where $C^\infty (M)'$ resp. $C^\infty_0 (M)'$ denotes the topological 
duals with respect to the natural locally convex topologies of the smooth 
functions resp. of the smooth functions with compact support. 
Which one of these topological duals one could or should choose in 
general depends much on the domain of definition of the 
positive linear functional (which need not to be the whole star product 
algebra $\CNP{C^\infty (M)}$ as e.~g. in section \ref {ExampleIISec}). 
Hence we shall concentrate in the following mainly on the purely 
algebraic properties and postpone these functional analytical 
questions as an interesting open problem to some future work.

This scheme opens up a lot of possibilities to study candidates for
positive linear functionals from a purely geometrical point of view:
one can always start with a classical measure $\omega_0$ having support
in an interesting submanifold of the classical phase space. The simplest
possibility consists in studying single points, that is the evaluation
or {\em Dirac delta functionals}, which we shall do in
section \ref {WickSec} and \ref {KaehlerSec}. 
One might also think about larger submanifolds such as certain
energy surfaces of suitable classical Hamiltonian functions, or Lagrangean
submanifolds. However, as we had already seen in section \ref {MotivSec} 
the positivity
of these functionals for the deformed algebra is in general no longer true,
and higher order terms in the deformation parameter have to be added.

One might also think about a deformed version of statistical mechanics
where $\omega_0$ is a classical statistical measure such as e.g. the
Boltzmann measure $exp(-\beta H_0)$ for a Hamiltonian function $H_0$ on
the phase space such that the Liouville integral over this Boltzmann 
function converges. In a context where convergence in $\lambda$ in a 
distributional sense is (at least partially) assumed one can find some 
of these ideas already in \cite{BFLS84}, \cite {BL85} where among other 
things KMS conditions are studied.
In the more particular case of flat $\mathbb R^{2n}$ and the 
convergent setting as mentioned in the introduction  
Hansen has proved convergence of star exponentials 
(Eq. (5.2) in \cite {Han84}) (of the above Boltzmann form) of tempered 
distributions satisfying a certain boundedness condition 
(see Eq. (5.5) and Corollary 5.2 in \cite {Han84}).
Note that in that paper GNS representations with vanishing Gel'fand ideal 
are considered, see Proposition 3.1.

\section {A remark on the time development in the GNS (pre-)Hilbert
          space}
\label {TimeSec}

In this section we should like to discuss the time development
in the observable algebra and in a representation generated by
some Hamiltonian. All results in this section are also valid 
(after the obvious modifications) when CNP-series are replaced 
by formal Laurent series.

Let $H_0 \in C^\infty (M)$ be a classical real Hamiltonian and let
$h \in \CNP{C^\infty (M)}$ be a real CNP series with $o(h) > 0$ and
define $H := H_0 + h$. Then the {\em Heisenberg equation} of motion with
respect to the Hamiltonian $H$ is given by
\BEQ {Heisenberg}
    \frac{d}{dt} f_t = \frac{1}{i\lambda} [f_t, H]
                     = \frac{1}{i\lambda} (f_t * H - H * f_t)
\EEQ
where $t \mapsto f_t \in \CNP{C^\infty (M)}$ is the trajectory through
$f_0 \in \CNP{C^\infty (M)}$ and $t \in \mathbb R$. We ask now for
solutions of (\ref {Heisenberg}) for a given initial condition $f_0$.
Let $X_{H_0}$ be the Hamiltonian vector field of $H_0$ and assume for
simplicity that $X_{H_0}$ is complete. In this case the flow
$\phi_t : M \to M$ of $X_{H_0}$ is a one-parameter group of symplectic
diffeomorphisms of M.
\begin {proposition}
With the above notations we have: The Heisenberg equation of motion has
a unique solution $f_t$ defined for all $t \in \mathbb R$ for any given
initial condition $f_0 \in \CNP{C^\infty (M)}$ and $f_t$ is given by
\BEQ {ZeitEntw}
    f_t = \phi^*_t \circ \mathcal T \exp \left(
    \int_0^t \phi^*_{-\tau} \circ \widehat H \circ \phi^*_\tau d\tau
    \right) f_0
\EEQ
where $\mathcal T\exp$ means the time ordered exponential and
$\widehat H (g) := \frac{1}{i\lambda} [g,H] - \{g,H_0\}$ satisfies
$o(\widehat H(g)) > o(g)$ for all $g \in \CNP{C^\infty (M)}$.
\end {proposition}
Remarks: Fedosov has proved a similar result in
\cite [Sec. 5.4.] {Fed96}. In the case where $X_{H_0}$ is not complete
one can prove the existence and uniqueness of solutions of
(\ref {Heisenberg}) for example for initial conditions $f_0$ with
compact support and sufficiently small times. We will denote the
solutions of (\ref{Heisenberg}) in the symbolical way
\[
    f_t = e^{\frac{it}{\lambda} \ad(H)} f_0
\]
where $\ad(H)f := [H, f]$.

We shall now turn to the time development of states in a specific
GNS \mbox{(pre-)}Hilbert space: We consider a positive linear functional
$\omega$ of $\CNP{C^\infty (M)}$ and the related GNS representation
$\pi_\omega$ on
$\mathfrak H_\omega := \CNP{C^\infty (M)}/\mathcal J_\omega$.
Then the Schr\"odinger equation is given by
\BEQ {Schroedinger}
    i\lambda \frac{d}{dt} \psi(t) = \pi_\omega (H) \psi (t)
\EEQ
where $\psi(t) \in \mathfrak H_\omega$. Specifying an initial condition
$\psi(0)$ one can rewrite the Schr\"odinger equation into an integral
equation
\[
    \psi (t) = \psi (0) + \frac{1}{i\lambda} \int_0^t
    \pi_\omega (H) \psi(\tau) d\tau
\]
but now iterating this equation will generate in general arbitrarily
high negative powers of $\lambda$ and hence it will not converge in the
formal sense! A rather simple example is provided by taking $H$ equal to
the constant function $1$.

But under certain conditions the well-defined time development of the
observables can be used to find solutions of (\ref{Schroedinger}):
\begin {theorem}
Let $\omega: \CNP{C^\infty (M)} \to \CNPC$ be a positive linear
functional and $H = H_0 + h \in \CNP{C^\infty (M)}$ as above. If
$H$ is contained in the Gel'fand ideal $\mathcal J_\omega$ of $\omega$
then the Schr\"odinger equation
(\ref {Schroedinger}) has a unique solution for $t \in \mathbb R$ and
any initial condition $\psi(0) \in \mathfrak H_\omega$. Moreover if
$f \in \CNP{C^\infty (M)}$ and $f_t$ is the solution of
(\ref {Heisenberg}) with initial condition $f$ then
\[
    \psi (t) := \psi_{f_{-t}} =
    \psi_{\exp(-\frac{it}{\lambda} \ad (H)) f}
\]
is the unique solution with initial condition $\psi (0) = \psi_f$.
\end {theorem}
\begin {proof}
Let $f_t$ be the solution of (\ref {Heisenberg}) with initial condition
$f$. Since $H \in \mathcal J_\omega$ we have
\[  
    \frac{d}{dt} \psi(t)
    = \frac{d}{dt} f_{-t} \bmod \mathcal J_\omega
    = - \frac{1}{i\lambda} (f_{-t} * H - H * f_{-t})
      \bmod \mathcal J_\omega
    = \frac{1}{i\lambda} \pi_\omega (H) \psi (t) .
\]
\end {proof}

Remarks: With the GNS representation we can define a time development of
the states in a purely formal way without considering any convergence
properties for $\lambda = \hbar$.
There are other methods to define a time development
as for example the {\em star exponential} \cite {BFFLS78} which is the
starting point for a spectral analysis. In this approach the time
development operator is viewed as a certain distribution if
one substitutes the formal parameter $\lambda$
by $\hbar \in \mathbb R$.

If $\omega$ is a state then the condition $H \in \mathcal J_\omega$
implies that the vacuum vector $\psi_{\bfmath 1}$ is invariant under
the time development and hence $\omega$ could be called an
{\em invariant state under $H$}. The search for invariant states under
certain group actions is a well-studied problem in algebraic quantum
field theory \cite {Haa93}. In our approach it raises the question
whether there is a positive linear functional such that the left ideal
generated by some Hamiltonian is contained in the Gel'fand ideal of that
functional. This condition would determine a preferred choice of a
state if one wants to consider the quantum theory of a certain classical
Hamiltonian.

\section {Example I:
          The Wick product and $\delta$-functionals in $\mathbb C^n$}
\label {WickSec}

We consider the Wick product (see for example \cite {BBEW96})
in $\mathbb C^n$
which is defined for $f,g \in \CNP{C^\infty (\mathbb C^n)}$ by
\BEQ {WickProdDef}
    f * g := \sum_{r=0}^\infty \frac{(2\lambda)^r}{r!}
             \frac {\partial^r f}
             {\partial z^{i_1} \cdots \partial z^{i_r}}
             \frac {\partial^r g}
             {\partial \cc z^{i_1} \cdots \partial \cc z^{i_r}}
\EEQ
where $z^1, \ldots, z^n$ are the canonical holomorphic coordinates in
$\mathbb C^n$ and summation over repeated indices is always understood.
Moreover we consider the evaluation functional (`delta-functional')
$\delta_p$ at the point $p \in \mathbb C^n$
\BEQ {DeltaDef}
    \delta_p [f] := f(p)
\EEQ
as a $\CNPC$-linear functional
$\delta_p : \CNP{C^\infty (\mathbb C^n)} \to \CNPC$. Then $\delta_p$
turns out to be positive with respect to the Wick product:
\begin {lemma} \label {DeltaWickPosLem}
Let $p \in \mathbb C^n$. Then $\delta_p$ is a positive linear functional
with respect to the Wick product (\ref {WickProdDef}). Moreover it is
clearly a state. The Gel'fand ideal of $\delta_p$ is given by
\BEQ {DeltaGelfand}
    \mathcal J_p := \left\{ f \in \CNP{C^\infty(\mathbb C^n)}
                    \; \left| \;
                    \frac {\partial^{|I|} f}{\partial \cc z^I} (p) = 0
                    \mbox { for all multiindices } I \right\} \right. .
\EEQ
\end {lemma}
This lemma is proved by observing that $\cc f * f$ evaluated at a point
is just a series of non-negative elements in $\CNPR$.
Using Borel's lemma \cite {Whi34} we easily find the
following isomorphism:
\begin {proposition}
We have the following isomorphism of $\CNPC$-vector spaces
\BEQ {Hp}
    \mathfrak H_p := \CNP{C^\infty (\mathbb C^n)} / \mathcal J_p
    \cong
    \CNP{ \left(\mathbb C[[\cc y^1, \ldots, \cc y^n]]\right)}
\EEQ
and an isomorhism is given by the formal $\cc z$-Taylor series at $p$
\BEQ {FormalTaylor}
    \mathfrak H_p \ni \psi_f \mapsto
    \sum_{r=0}^\infty \frac {1}{r!}
    \frac{\partial^r f}
    {\partial \cc z^{i_1} \cdots \partial \cc z^{i_r}} (p)
    \cc y^{i_1} \cdots \cc y^{i_r}
\EEQ
and the Hermitian product induced by $\delta_p$ is given by
\BEQ {HermProd}
    \SP {\psi_f} {\psi_g} 
    = \sum_{r=0}^\infty \frac{(2\lambda)^r}{r!}
      \cc {\frac{\partial^r f}
      {\partial \cc z^{i_1} \cdots \partial \cc z^{i_r}} (p)}
      \frac{\partial^r g}
      {\partial \cc z^{i_1} \cdots \cc z^{i_r}}(p) .
\EEQ
\end {proposition}
In the following we shall identify $\psi_f$ with the formal $\cc
z$-Taylor series and consider only $p=0$ for simplicity. Then we can
determine the GNS representation $\pi_0$ induced by $\delta_0$ by
calculating the formal $\cc z$-Taylor series of $f*g$ in order to obtain
$\pi_0(f)\psi_g$.
\begin {lemma} [The formal Wick representation]
For $f \in \CNP{C^\infty (\mathbb C^n)}$ we find
\BEQ {pinullf}
    \pi_0 (f) = \sum_{r,s=0}^\infty \frac{(2\lambda)^r}{r!s!}
    \frac{\partial^{r+s} f}{\partial z^{i_1} \cdots \partial z^{i_r}
    \partial \cc z^{j_1} \cdots \partial \cc z^{j_s}} (0)
    \cc y^{j_1} \cdots \cc y^{j_s}
    \frac{\partial}{\partial \cc y^{i_1}} \cdots
    \frac{\partial}{\partial \cc y^{i_r}}
\EEQ
In particular we have for polynomials the Wick ordering (normal ordering):
\BEQ {WickPoly}
    \pi_0 (z^{i_1} \cdots z^{i_r} \cc z^{j_1} \cdots \cc z^{j_s})
    =
    (2\lambda)^r \cc y^{j_1} \cdots \cc y^{j_s}
    \frac{\partial}{\partial \cc y^{i_1}} \cdots
    \frac{\partial}{\partial \cc y^{i_r}}
\EEQ
\end {lemma}

Now we will consider the completion of $\mathfrak H_0$ to a
$\CNPC$-Hilbert space.
\begin {proposition} \label {WickHProp}
The completion $\widehat{\mathfrak H_0}$ of $\mathfrak H_0$
is given by
\begin {eqnarray}
    \widehat{\mathfrak H_0}
    & = & \left\{ \phi = \sum_{r=0}^\infty
          \frac{1}{r!} a_{i_1 \cdots i_r} (\lambda)
          \cc y^{i_1} \cdots \cc y^{i_r}
          \; \right| \;
          a_{i_1 \cdots i_r} (\lambda) \in \CNPC \nonumber \\
    &   & \label {CompHilbert}
          \mbox { such that } 
          \left. \sum_{r=0}^\infty \frac{(2\lambda)^r}{r!}
          \cc {a_{i_1\cdots i_r} (\lambda)} a_{i_1\cdots i_r} (\lambda)
          \mbox { converges in } \CNPC \right\} .
          \label {CompleteDef}
\end {eqnarray}
Let $K = (k_1, \ldots, k_n)$ be a multiindex then the vectors
\BEQ {HilbertBase}
    e_K := \frac{1}{\sqrt{\lambda^{|K|} K!}}
           (\cc y^1)^{k_1} \cdots (\cc y^n)^{k_n}
    \in \widehat{\mathfrak H_0}
    \qquad \mbox { where }
    \begin {array} {c}
    |K| := k_1 + \cdots + k_n \\
    K! := k_1! \cdots k_n!
    \end {array}
\EEQ
form a {\em Hilbert base} for $\widehat{\mathfrak H_0}$ and hence
$\widehat{\mathfrak H_0}$ is isometric to $\ell^2(\CNPC)$.
\end {proposition}
Remark: An element in $\mathfrak H_0$ has a lowest order in $\lambda$
according to equation (\ref {Hp}) but the coefficients
$a_{i_1 \cdots i_r} (\lambda)$ of an element in
$\widehat{\mathfrak H_0}$ written as in (\ref {CompleteDef})
could have decreasing orders in $\lambda$. Assigning to
$\cc y^1, \ldots, \cc y^n$ the total degree $1$ and to $\lambda$ the
total degree $2$ the elements in $\widehat{\mathfrak H_0}$ could be
understood as formal CNP series in this total degree.

As an example we want to consider the harmonic oscillator with
classical Hamiltonian $H:= \frac{1}{2} \omega |z|^2$ where 
$\omega \in \mathbb R^+$ is the oscillator frequency. First
we notice that $H$ is an element of the Gel'fand ideal $\mathcal J_0$ 
and $\pi_0(H) = \lambda \omega \cc y^k \partial/\partial \cc y^k$.
Note that this operator is clearly defined on all of
$\widehat{\mathfrak H_0}$ as a bounded and hence continuous operator.
We can ask for its {\em GNS-spectrum} in the following sense:
\[
    \spec(\pi_0 (H)) := \{\mu \in \CNPC \; | \;
    (\pi_0(H) - \mu \bfmath 1) \mbox{ is bijective } \}
\]
By an easy calculation we find that the GNS-spectrum
is purely discrete, namely:
\[
    \spec(\pi_0(H)) = \{0, \lambda\omega, 2\lambda\omega, \ldots \}
\]
Moreover we notice that the vectors $e_K$ are eigenvectors 
to the eigenvalue $\lambda\omega|K|$ and hence they are a Hilbert base of 
eigenvectors to the harmonic oscillator.
Note that the GNS-spectrum is still formal (i. e. a `CNP-number')
by definition, but after the substitution $\lambda = \hbar$ we arrive at
the well-known real oscillator spectrum (minus the ground state energy
which is due to the Wick product) which --in the context of star
products-- had directly been computed by analytic methods using
the star exponential in \cite [p. 125--132] {BFFLS78}.

At last we want to discuss a way how one can get back convergence in
$\lambda$ if we substitute $\lambda$ by $\hbar \in \mathbb R$. The main
idea is that we ask for convergence in the representation and {\em not}
in the algebra. We define for a fixed $\hbar \in \mathbb R^+$
\begin {eqnarray}
    H(\hbar)
    & := & \left\{ \phi \in \widehat{\mathfrak H_0} \; \left | \;
           \sum_{r=0}^\infty \hbar^{q_r} a_{q_r}^{(K)} 
           := \SP {e_K} \phi \big|_{\lambda = \hbar}
           \mbox { converges } \right.\right.  \\
    &    & \left. \mbox { absolutely } \forall K \mbox{ and }
           \sum_{|K|=0}^\infty
           \left(
           \sum_{r=0}^\infty \hbar^{q_r} \left|a_{q_r}^{(K)}\right| \right)^2
           < \infty \right\} \label {HhbarDef} \\
    N(\hbar)
    & := & \left\{\left. \phi \in \widehat{\mathfrak H_0} \; \right| \;
           \SP{e_K} \phi \big|_{\lambda = \hbar} = 0 
           \right\} .
\end {eqnarray}
\begin {lemma} \label {ConvLem}
Let $\phi = \sum_{r=0}^\infty \frac{1}{r!} a_{i_1 \ldots i_r} (\lambda)
\cc y^{i_1} \cdots \cc y^{i_r} \in \widehat{\mathfrak H_0}$.
If $\phi \in H(\hbar)$ then we have:
\begin {enumerate}
\item $a_{i_1 \ldots i_r}(\lambda)$ converges absolutely
      for $\lambda = \hbar$.
\item The Hermitian product $\SP \phi\phi$ is absolutely convergent for
      $\lambda = \hbar$ and positive semidefinite in $\mathbb C$:
      \[
          \SP \phi\phi \big|_{\lambda = \hbar} =
          \sum_{r=0}^\infty \frac{(2\hbar)^r}{r!}
          \cc {a_{i_1 \ldots i_r} (\hbar)}
          a_{i_1 \ldots i_r} (\hbar)
          \ge 0
      \]
\item $H(\hbar)$ is a $\mathbb C$-vector space and $N(\hbar)$ is a
      $\mathbb C$-subvector space of $H(\hbar)$. For
      $\phi, \psi \in H(\hbar)$ the Hermitian product $\SP \phi \psi$
      converges absolutely in $\mathbb C$ for $\lambda = \hbar$.
\item $F(z) := \sum_{r=0}^\infty \frac{a_{i_1 \ldots i_r} (\hbar)}{r!}
      \cc z^{i_1} \cdots \cc z^{i_r}$ is an entire
      anti-holomorphic function on $\mathbb C^n$.
\end {enumerate}
\end {lemma}
Part {\it ii.)} and {\it iii.)} will be proved in a more
general context in theorem \ref {HNConvTheo}.
Let $\phi$ be given as in the above lemma and let
$\psi = \sum_{r=0}^\infty \frac{1}{r!} b_{i_1 \ldots i_r} (\lambda)
\cc y^{i_1} \cdots \cc y^{i_r}$ be another element in
$H(\hbar)$. Then we have the following theorem:
\begin {theorem}
The quotient $\mathbb C$-vector space
$\mathcal H(\hbar) := H(\hbar)/N(\hbar)$ with the Hermitian product
\BEQ {ConvProdDef}
    \SP {\phi \bmod N(\hbar)} {\psi \bmod N(\hbar)} :=
    \sum_{r=0}^\infty \frac{(2\hbar)^r}{r!}
    \cc {a_{i_1 \ldots i_r} (\hbar)} b_{i_1 \ldots i_r} (\hbar)
\EEQ
is a Hilbert space over $\mathbb C$ which is isometric to the Hilbert
space of {\em anti-holomorphic} functions $f:\mathbb C^n \to \mathbb C$
\BEQ {BargmannHilbDef}
    \mathcal H(\hbar) \cong
    \left\{ f: \mathbb C^n \to \mathbb C
    \; \left| \; \sum_{r=0}^\infty \frac{(2\hbar)^r}{r!}
    \sum_{i_1,\ldots,i_r} \left|\frac{\partial^r f}
    {\partial \cc z^{i_1} \cdots \partial \cc z^{i_r}} (0) \right|^2
    < \infty
    \right\} \right.
\EEQ
with the Gaussian product
\BEQ {GausprodDef}
    \SP f g =
    \frac{1}{(2\pi\hbar)^n} \int_{\mathbb C^n} e^{-\frac{|z|^2}{2\hbar}}
    \cc{f(z)} g(z)
    dz^1 \cdots dz^n d\cc z^1 \cdots d\cc z^n,
\EEQ
and the canonical isometry is given by $\phi \mapsto F$ where $F$ is the
anti-holomorphic function defined by $\phi$ as in
lemma \ref {ConvLem} iv.).
\end {theorem}

Moreover the GNS representation of $\CNP{C^\infty (\mathbb C^n)}$
induces a representation $\mathcal Q$ on the $\mathbb C$-Hilbert space
$\mathcal H(\hbar)$ at least for `many' functions. We consider only
those elements $f \in \CNP{C^\infty (\mathbb C^n)}$ with
$\pi_0(f) N(\hbar) \subseteq N(\hbar)$. Then we define
\BEQ {DfDef}
    D_f(\hbar) := \{\phi \in H(\hbar) \; | \;
                  \pi_0(f) \phi \in H(\hbar) \} \; ;
    \qquad
    \mathcal D_f (\hbar) := D_f (\hbar) / N(\hbar)
\EEQ
and call $\mathcal D_f (\hbar) \subseteq \mathcal H(\hbar)$ the domain of
$f$ since we can define an operator $\mathcal Q (f)$ with domain of
definition $\mathcal D_f (\hbar)$ by
\BEQ {WickQuantDef}
    \mathcal Q(f) (\phi \bmod N(\hbar)) := \pi_0 (f) \phi \bmod N(\hbar)
    \qquad
    \mbox { for } \phi \in D_f (\hbar) .
\EEQ
In general $\mathcal D_f(\hbar)$ may be very small
but there are `many' functions (i.e. all polynomials in $z$ and $\cc z$) 
such
that $\mathcal D_f(\hbar)$ is dense in $\mathcal H(\hbar)$ or even equal
to $\mathcal H(\hbar)$. In this case we clearly have
\BEQ {WickRep}
    \mathcal Q (f*g) = \mathcal Q (f) \circ \mathcal Q (g)
\EEQ
and hence we call $\mathcal Q$ a {\em quantization map}. If we write
$\mathcal H(\hbar)$ in the form (\ref{BargmannHilbDef}) this is of
course just the well-known Bargmann quantization.

\section {General results for K\"ahler manifolds}
\label {KaehlerSec}

In this section we shall derive some general results for K\"ahler
manifolds with the Fedosov star product of Wick type which had been
constructed in \cite {BW96a}.
First we define the {\em formal Wick algebra} in $2n$ parameters by
\BEQ {FormWickDef}
    \mathcal W_n := \CNP{\left( \mathbb C[[y^1, \ldots, y^n,
                    \cc y^1, \ldots, \cc y^n]] \right)}
\EEQ
then $\mathcal W_n$ is clearly a $\CNPC$-vector space and we
define a pointwise multiplication of the formal power series in
$y^1, \ldots, \cc y^n$ as usual. Then $\mathcal W_n$ becomes an
associative and commutative algebra over $\CNPC$ with unit element. 
Moreover $\mathcal W_n$ becomes a $^*$-algebra if we extend the complex
conjugation such that $y^k$ is mapped to $\cc y^k$ and vice versa.
This formal Wick algebra could be deformed by an analogue of the usual
Wick product. We define
\BEQ {FormWickProdDef}
    a \circ b := \sum_{r=0}^\infty \frac{(2\lambda)^r}{r!}
                 \frac{\partial^r a}
                 {\partial y^{i_1} \cdots \partial y^{i_r}}
                 \frac{\partial^r b}
                 {\partial \cc y^{i_1} \cdots \partial \cc y^{i_r}}
\EEQ
for $a,b \in \mathcal W_n$. This leads to an associative deformation of
the pointwise multiplication such that the complex conjugation is still
an antilinear algebra involution and $a \circ 1 = 1 \circ a = a$.
Every element $a \in \mathcal W_n$ could be written as
\[
    a = \sum_{r,s=0}^\infty \frac{1}{r!s!}
        a_{i_1 \ldots i_r j_1 \ldots j_s}
        y^{i_1} \cdots y^{i_r} \cc y^{j_1} \cdots \cc y^{j_s} .
\]
Then we define the delta functional as the projection onto
the part without explicit powers of the formal parameters
$y^1, \ldots, \cc y^n$:
\[
    \delta (a) := a_{00} \in \CNPC
\]
This will also be written in the following symbolical way:
$\delta(a) = a_{00} = a|_{y=0}$. Then it is easy to see that the delta
functional is a positive linear functional:
\begin {proposition} \label {FormWickRepProp}
\begin {enumerate}
\item The delta functional $\delta: \mathcal W_n \to \CNPC$ is
      a state of $\mathcal W_n$ with respect to the Wick product
      and the Gel'fand ideal of $\delta$ is given by:
      \BEQ {WickGel}
          \mathcal J = \left\{ a \in \mathcal W_n \; \left| \;
                       \left.\frac{\partial^{|I|} a}{\partial \cc
                       y^I}\right|_{y=0} = 0
                       \mbox{ for all multiindices } I \right\} \right.
      \EEQ
\item The quotient space $\mathcal W_n/\mathcal J$ is isomorphic to
      $\mathfrak H'_n :=
      \CNP{\left(\mathbb C[[\cc y^1, \ldots, \cc y^n]] \right)}$
      with the `formal $\cc y$-Taylor series' as isomorphism
      \[
          a \bmod \mathcal J \mapsto
          \sum_{s=0}^\infty \frac{1}{s!}
          \left. \frac{\partial^s a}
          {\partial \cc y^{i_1} \cdots \partial \cc y^{i_s}}
          \right|_{y=0}
          \cc y^{i_r} \cdots \cc y^{i_s}
      \]
      and the Hermitian product induced by $\delta$ is given by
      \[
          \SP{\sum_{r=0}^\infty \frac{1}{r!}
              a_{i_1 \ldots i_r} \cc y^{i_r} \cdots \cc y^{i_r}}
             {\sum_{s=0}^\infty \frac{1}{s!}
             b_{j_1 \ldots j_s} \cc y^{j_r} \cdots \cc y^{j_s}}
          =
          \sum_{r=0}^\infty \frac{(2\lambda)^r}{r!}
          \cc{a_{{i_1}\ldots{i_r}}} b_{{i_1}\ldots{i_r}}
      \]
\end {enumerate}
\end {proposition}
The properties of $\mathfrak H'_n$ and its completion $\mathfrak H_n$ 
were already discussed in proposition \ref {WickHProp} where
$\mathfrak H'_n \cong \mathfrak H_0$.

Now we want to investigate the convergence of the
Hermitian products in a Hilbert space $\mathfrak H$ over $\CNPC$ with a
countable Hilbert base $\{e_k\}_{k \in \mathbb N}$. As we will see in
section \ref {ExampleIISec} a $\CNPC$-Hilbert space has in general no
Hilbert base but we will
see important examples if we consider K\"ahler manifolds.

First we will mention some easy convergence properties of CNP series. Let
$\mathbb C^- := \mathbb C \setminus \{x \in \mathbb R \; | \; x \le 0 \}$
and $B_r^- := B_r(0) \cap \mathbb C^-$ where $B_r(0)$ is the open disk
around $0 \in \mathbb C$ of radius $r$. For $q \in \mathbb Q$ let $z^q$
be the holomorphic root defined on $\mathbb C^-$. For a formal CNP
series $a = \sum_{r=0}^\infty \lambda^{q_r} a_{q_r} \in \CNPC$ we define
the radius of convergence by
$R := \sup_{t \in \mathbb R^+}
\{\sum_{r=0}^\infty t^{q_r} |a_{q_r}| < \infty \}$
then we can prove the following proposition:
\begin {proposition}
Let $a \in \CNPC$ with convergence radius $R > 0$.
Then $f(z) := \sum_{r=0}^\infty z^{q_r} a_{q_r}$
converges absolutely and normally in $B_R^-$. Hence
$f: z \in B_R^- \to f(z) \in \mathbb C$ is a holomorphic function.
\end {proposition}
Now we define $H(\hbar)$ and $N(\hbar)$ in the following way:
\begin {eqnarray}
    H(\hbar)
    & := & \left\{ \phi \in \mathfrak H \; \left| \;
           \forall k \in \mathbb N :
           \SP {e_k}\phi\big|_{\lambda = \hbar} = \sum_{r=0}^\infty
           \hbar^{q_r} a_{q_r}^{(k)}
           \mbox { converges absolutely } \right.\right. \nonumber \\
    &    & \left. \mbox { and }
           \sum_{r=0}^\infty \left| F_k (\hbar) \right|^2 < \infty
           \mbox { where }
           F_k (\hbar) := \sum_{r=0}^\infty \hbar^{q_r}
           \left|a_{q_r}^{(k)}\right| \right\} \\
    N(\hbar)
    & := & \left\{ \phi \in H(\hbar) \; \left| \;
           \forall k \in \mathbb N :
           \SP {e_k}\phi \big|_{\lambda = \hbar} = 0 \right\} \right.
\end {eqnarray}
First we notice that $F_k (\hbar)$ is well-defined and
extends to a holomorphic function $F_k$ on $B_\hbar^-$. Using
Cauchy's theorem for iterated series and Weierstra\ss' theorem
for convergence of holomorphic functions we easily get the
following theorem:
\begin {theorem} \label {HNConvTheo}
Let $\mathfrak H$ be a $\CNPC$-Hilbert space with countable
Hilbert base $\{ e_k \}_{k \in \mathbb N}$ and let $H(\hbar)$ and
$N(\hbar)$ be defined as above and $\phi, \psi \in H(\hbar)$.
\begin {enumerate}
\item $H(\hbar)$ is a $\mathbb C$-vector space and $N(\hbar)$ is a
      $\mathbb C$-subvector space of $H(\hbar)$.
\item The Hermitian product of $\mathfrak H$ induces a semidefinite
      sesquilinear form on $H(\hbar)$ by
      \BEQ {SPhbarDef}
          \SP \phi \psi \Big|_{\lambda = \hbar} =
          {\SP \phi \psi}_\hbar :=
          \sum_{k=1}^\infty \SP \phi {e_k} \big|_{\lambda = \hbar}
          \SP {e_k}\psi \big|_{\lambda = \hbar}
      \EEQ
      and ${\SP \phi\phi}_\hbar = 0$ iff $\phi \in N(\hbar)$.
\item $\SP \phi \psi \big|_{\lambda = z}$ is a holomorphic function on
      $B_\hbar^-$ and
      \[
          \SP \phi \psi \big|_{\lambda = z} = \sum_{k=1}^\infty
          \SP \phi {e_k} \big|_{\lambda = z}
          \SP {e_k}\psi \big|_{\lambda = z}
      \]
      converges normally on $B_\hbar^-$.
\item The quotient space $\mathcal H(\hbar) := H(\hbar) / N(\hbar)$ with
      the induced Hermitian product is a $\mathbb C$-Hilbert space with
      Hilbert base $\{e_k \bmod N(\hbar) \}_{k \in \mathbb N}$ and hence
      $\mathcal H(\hbar)$ is isometric to $\ell^2(\mathbb C)$.
\end {enumerate}
\end {theorem}
Note that the resulting $\mathbb C$-Hilbert space
$\mathcal H$ is already a Hilbert space and not only a pre-Hilbert
space.

Remark: If one defines $F_k$ with $a^{(k)}_{q_r}$ instead of its
absolute value then one would obtain again two $\mathbb C$-subvector
spaces and again the quotient would be a Hilbert space. But now
$\SP \phi\psi |_{\lambda = \hbar}$ is in general not convergent to
${\SP \phi\psi}_\hbar$ as in equation (\ref {SPhbarDef}). Consider for
example
$\phi_1 := \sum_{k=0}^\infty \lambda^k \sin(\frac{\pi \lambda}{2}) e_k$
or $\phi_2 := \sum_{r=0}^\infty k!
\lambda^k \sin(\frac{\pi \lambda}{2}) e_k$ which are both in $H(2)$ but
clearly $\SP {\phi_1}{\phi_1}|_{\lambda = \hbar}$ has only convergence
radius $1$ and $\SP {\phi_2}{\phi_2}|_{\lambda = \hbar}$ has convergence
radius $0$ while
${\SP {\phi_1}{\phi_1}}_\hbar = {\SP {\phi_2} {\phi_2}}_\hbar = 0$.

Now we want to consider the Fedosov star product of Wick type for a
K\"ahler manifold $M$ of real dimension $2n$ constructed as in
\cite {BW96a}. We will use the same notation as in \cite {BW96a} with the
only difference that the formal parameter in denoted by $\lambda$ and
$\hbar$ is reserved for the real number corresponding to the value of
Planck's constant in a chosen unit system.
Let $\omega$ be the symplectic K\"ahler form on $M$ which
is given in local holomorphic coordinates by
$\omega = \frac{i}{2} \omega_{k \cc l} dz^k \wedge d\cc z^l$.
The Fedosov algebra is defined by
\BEQ {FedAlgDef}
    \WL := \CNP{\left({\mathsf X}_{s=0}^\infty \mathbb C\left(
           \Gamma\left(\mbox{$\bigvee$}^s T^*M
           \otimes \mbox{$\bigwedge$} T^*M\right)\right)\right)}
\EEQ
together with the pointwise multiplication induced by the symmetric and
antisymmetric product of forms.
The elements $a \in \WL$ are of the form
$\sum_{r=0}^\infty \lambda^{q_r} a_{q_r}$ where
$a_{q_r} = \sum_{s=0}^\infty a^{(q_r)}_s$ and
$a^{(q_r)}_s \in
\mathbb C (\Gamma (\bigvee^s T^*M \otimes \bigwedge T^*M))$ are smooth
sections. The fibrewise Wick product as a deformation of the pointwise
product is defined by
\BEQ {FibWickDef}
    \begin {array} {c}
    \displaystyle
    a \circ b := \sum_{r=0}^\infty \left(\frac{i\lambda}{2}\right)^r
                 {\Lambda}^{(r)} (a,b) \\
    \displaystyle
    {\Lambda}^{(r)} (a,b) :=
                 \frac{1}{r!} \left(\frac{4}{i}\right)^r
                 \omega^{k_1\cc l_1} \cdots \omega^{k_r\cc l_r}
                 i_s(Z_{k_1}) \cdots i_s (Z_{k_r}) a \;
                 i_s(\cc Z_{l_1}) \cdots i_s (\cc Z_{l_r}) b .
    \end {array}
\EEQ
and since this product is defined fibrewise we can define the Fedosov
algebra at a point $p \in M$ by
\BEQ {FedAlgpDef}
    \WL_p := \CNP{\left({\mathsf X}_{s=0}^\infty
             \mathbb C\left(\mbox{$\bigvee$}^s T^*_pM
             \otimes \mbox{$\bigwedge$} T^*_pM\right)\right)}
\EEQ
together with the product $\circ$. Then the restriction $a_p$
of a section $a \in \WL$ to the point $p \in M$ is an element in
$\WL_p$. The sections in $\WL$ without any antisymmetric part are
denoted by $\W$ and analogously we define $\W_p$.
In \cite {BW96a} we have shown that there exists a
{\em Fedosov derivation} $D$ for $\WL$ such that $D^2 = 0$
and $\W_D := \ker D \cap \W$ is
isomorphic to $\CNP{C^\infty (M)}$ and the isomorphism is given by
the {\em Fedosov-Taylor series} $\tau: \CNP{C^\infty (M)} \to \W_D$
which was constructed recursively and the inverse map
$\sigma: \W_D \to \CNP{C^\infty (M)}$ is simply the projection on the
elements with symmetric degree zero, i.e. $\sigma = \pi^{(0,0)}$ where
$\pi^{(r,s)}$ is the projection onto the symmetric forms of type $(r,s)$.
Then the Fedosov star product is given by
$f * g = \sigma (\tau(f) \circ \tau (g))$ and several properties of this
star product were shown in \cite {BW96a}. The most important property for
our purpose is the {\em reality} of the Fedosov-Taylor series
$\tau(\cc f) = \cc {\tau (f)}$.
\begin {lemma} \label {IsoWpWn}
Let $M$ be a $2n$-dimensional K\"ahler manifold and
$z^1, \ldots, z^n$ a
holomorphic chart around $p \in M$ such that
$\omega_{k \cc l}\big|_p = \delta_{k \cc l}$ in this chart. Then
$\varphi (1) := 1$, $\varphi (dz^k) := y^k$,
$\varphi (d\cc z^k) := \cc y^k$ induces a $\CNPC$-algebra
isomorphism $\varphi: \W_p \to \W_n$ with
respect to the fibrewise Wick product in $\W_p$ and the
Wick product (\ref {FormWickProdDef}) in $\W_n$.
\end {lemma}
Note that it is always possible to find such a holomorphic chart
\cite [Sec. 0.7] {GH78}. With such an isomorphism we easily find the
following proposition analogously to lemma \ref {DeltaWickPosLem}:
\begin {proposition}
Let $p \in M$ and $\delta_p$ the delta functional at $p$.
Then $\delta_p: \CNP{C^\infty (M)} \to \CNPC$ is a state with respect to
the Fedosov star product of Wick type and the Gel'fand ideal is given by
\[
    \mathcal J_p := \left\{ f \in \CNP{C^\infty (M)} \; \Big| \;
                    \forall r\ge 0 : \pi^{(0,r)} \tau_p (f) = 0
                    \right\}
\]
where $\tau_p (f)$ is the Fedosov-Taylor series of $f$ evaluated at $p$.
\end {proposition}
In the next proposition we shall prove that the Fedosov-Taylor series
$\tau_p$ at a point $p \in M$ is surjective on $\W_p$ where we use
Borel's lemma and the recursion
formula \cite [eqn. (19)] {BW96a} for $\tau$.
\begin {proposition} \label {FedTaySurProp}
The Fedosov-Taylor series $\tau_p : \CNP{C^\infty (M)} \to \W_p$
is surjective for all $p \in M$.
\end {proposition}
We define for $p \in M$
\[
    \widetilde{\mathcal J}_p
    := \left\{ a \in \W_p \; \Big| \;
       \forall r \ge 0 :
       \pi^{(0,r)} a = 0 \right\}
    \qquad
    \widetilde{\mathfrak H'}_p
    := \CNP{\left({\mathsf X}_{s=0}^\infty
       \mbox{$\bigvee$}^{(0,s)} T^*_pM
       \right)}
\]
and then we can describe the GNS representation induced by $\delta_p$ in
the following way:
\begin {theorem}
With the notations from above we have:
\begin {enumerate}
\item The quotient space
      $\mathfrak H'_p := \CNP{C^\infty (M)}/\mathcal J_p$ is canonically
      isomorphic to $\widetilde{\mathfrak H'}_p$ where the isomorphism
      is induced by $\tau_p$ (and also denoted by $\tau_p$):
      \begin {eqnarray*}
          \psi_f \mapsto \tau_p (\psi_f)
          & := & \tau_p (f) \bmod \widetilde{\mathcal J}_p \\
          & =  & \sum_{r=0}^\infty \frac{1}{r!}
                 \left( \pi^{(0,0)} i_s (\cc Z_{k_1}) \cdots i_s (\cc Z_{k_r})
                 \tau_p (f) \right) d\cc z^{k_1} \vee \cdots \vee d\cc z^{k_r}
      \end {eqnarray*}
      where $\psi_f \in \mathfrak H'_p$
      and the Hermitian product is given by
      \begin {eqnarray*}
          \SP{\psi_f}{\psi_g}
          & = & \sum_{r=0}^\infty \frac{(2\lambda)^r}{r!}
                \omega^{k_1 \cc l_1}_p \cdots \omega^{k_r \cc l_r}_p \\
          &   & \times \; \pi^{(0,0)}
                \left(i_s (Z_{k_1}) \cdots i_s (Z_{k_r})
                \cc {\tau_p (f)} i_s (\cc Z_{l_1}) \cdots i_s (\cc Z_{l_r})
                \tau_p (g) \right)
      \end {eqnarray*}
      and the GNS representation induced by $\delta_p$ is given by
      \begin {eqnarray*}
          \lefteqn{\tau_p \circ \pi_p (f) \circ \tau_p^{-1} } \\
          & = & \sum_{r,s=0}^\infty \frac{(2\lambda)^r}{r!s!}
                \left(\pi^{(0,0)} \left(
                i_s (Z_{k_1}) \cdots i_s (Z_{k_r}) i_s (\cc Z_{l_1})
                \cdots i_s (Z_{l_s}) \tau_p (f)\right)\right) \\
          &   & \times \;
                \omega^{k_1 \cc t_1} \cdots \omega^{k_r \cc t_r}
                d\cc z^{l_1} \vee \cdots \vee d\cc z^{l_s}
                i_s (\cc Z_{t_1}) \cdots i_s (\cc Z_{t_r})
      \end {eqnarray*}
\item Let $z^1, \ldots, z^n$ be a holomorphic chart such that
      $\omega_{k\cc l}|_p = \delta_{k \cc l}$ and let
      $\varphi: \W_p \to \W_n$ be defined as in lemma \ref {IsoWpWn}.
      Then $\varphi \circ \tau_p$ induces an
      isomorphism (also denoted by $\varphi \circ \tau_p$)
      of the $\CNPC$-pre-Hilbert space $\mathfrak H'_p$ to
      $\mathfrak H'_n$ and hence the completion $\mathfrak H_p$
      of $\mathfrak H'_p$ is isometric to $\mathfrak H_n$
      via $\varphi \circ \tau_p$ and hence $\mathfrak H_p$ has a
      countable Hilbert base and is isometric to $\ell^2(\CNPC)$.
      Moreover we have for $\psi_f \in \mathfrak H'_p$
      \[
         \varphi \circ \tau_p (\psi_f) =
         \sum_{r=0}^\infty \frac {1}{r!}
         \left.\frac{\partial^r \varphi \circ \tau_p (f)}
         {\partial \cc y^{k_1} \cdots \partial \cc y^{k_r}}\right|_{y=0}
         \cc y^{k_1} \cdots \cc y^{k_r} .
      \]
\end {enumerate}
\end {theorem}

As in the Wick case we can ask for convergence. First we have to choose
a Hilbert base of $\mathfrak H_p$. In a holomorphic chart
$z^1, \ldots, z^n$ with $\omega_{k\cc l}|_p = \delta_{k\cc l}$ we can
use for example the vectors
$\hat e_K := (\varphi \circ \tau_p)^{-1} e_K$ with $e_K$ as in
proposition \ref {FormWickRepProp}. Then we construct $H_p(\hbar)$ and
$N_p(\hbar)$ with respect to this Hilbert base
as in theorem \ref {HNConvTheo} and get a $\mathbb C$-Hilbert space
$\mathcal H_p (\hbar)$ and a representation of `many' functions
$f \in \CNP{C^\infty (M)}$ in the same way as in the Wick case.
We consider again functions $f$ with
$\pi_p (f) N_p(\hbar) \subseteq N_p(\hbar)$ and define $D_f (\hbar)$ and
$\mathcal D_f(\hbar)$ as in the Wick case.
Then we define the quantization map $\mathcal Q$ analogously:
$\mathcal Q(f)$ is an operator on $\mathcal D_f (\hbar)$ defined by
\[
    \mathcal Q(f) (\psi \bmod N_p(\hbar))
    := \pi_p (f) \psi \bmod N_p(\hbar) \qquad \psi \in D_f (\hbar)
\]
which leads to the representation
property $\mathcal Q(f) \circ \mathcal Q(g) = \mathcal Q(f * g)$ on
suitable domains.
Since according to proposition \ref {FedTaySurProp} the Fedosov-Taylor
series is surjective there are indeed `many' suitable functions such as
the pre-images under $\tau_p$ of the polynomials in $dz^k$ and
$d\cc z^k$ (where $z^1, \ldots, z^n$ is the chart from above).

\section {Example II: The Weyl-Moyal product for $\mathbb R^{2n}$}
\label {ExampleIISec}

Now we want to consider the phase space
$T^*\mathbb R^n \cong \mathbb R^{2n}$ with the standard symplectic form
and the Weyl-Moyal product defined for
$f,g \in \CNP{C^\infty (\mathbb R^{2n})}$ by
\BEQ {WeylProdDef}
    \begin {array} {c}
    \displaystyle
    f * g := \sum_{r=0}^\infty \left(\frac{i\lambda}{2}\right)^r
             \Lambda^{(r)} (f,g) \\
    \displaystyle
    \Lambda^{(r)} (f,g) := \frac{1}{r!}
    \Lambda^{i_1j_1} \cdots \Lambda^{i_rj_r}
    \frac{\partial^r f}{\partial x^{i_1} \cdots \partial x^{i_r}}
    \frac{\partial^r g}{\partial x^{j_1} \cdots \partial x^{j_r}}
    \end {array}
\EEQ
where $x^1, \ldots, x^{2n} = q^1, \ldots, q^n, p_1, \ldots p_n$ and
$\Lambda^{ij}$ are the components of the Poisson tensor $\Lambda$ with
respect to the coordinates $x^1, \ldots, x^{2n}$, i.e.
$\Lambda = \frac{1}{2} \Lambda^{ij} \partial_{x^i} \wedge
\partial_{x^j} = \partial_{q^i} \wedge \partial_{p_i}$.
The smooth functions on $\mathbb R^{2n}$ with compact support are
denoted by $\CNP{C^\infty_0 (\mathbb R^{2n})}$.
\begin {lemma}
Let $f,g \in \CNP{C^\infty_0 (\mathbb R^{2n})}$ then the integral over
$\mathbb R^{2n}$ is a trace \cite {CFS92}
\BEQ {IntTrace}
    \int_{\mathbb R^{2n}} f*g \, d^{2n} x =
    \int_{\mathbb R^{2n}} g*f \, d^{2n} x =
    \int_{\mathbb R^{2n}} fg \, d^{2n} x
\EEQ
and clearly a positive linear functional.
\end {lemma}

To obtain Schr\"odinger's quantization we need a different positive
linear functional namely the 
{\em integration over the configuration space}
$Q := \mathbb R^n$ of $T^*Q \cong \mathbb R^{2n}$
for a fixed value $\vec{p_0}$ of the momenta. But first we have to
define a suitable subalgebra of
$\CNP{C^\infty (T^*Q)}$ such that the integration is
well-defined. For $\vec{p_0} \in \mathbb R^n$ we define
\BEQ {CpDef}
    C^\infty_{\vec{p_0}} (T^*Q) :=
    \left\{ f \in C^\infty (T^*Q) \; \left| \;
    \supp f \cap Q_{\vec{p_0}} \mbox { is compact } \right\} \right.
\EEQ
where
$Q_{\vec{p_0}} :=
\{ (\vec{q}, \vec{p_0}) \in T^*Q \; | \; \vec{q} \in Q \}$
and notice that due to the locality of the Weyl-Moyal product
$\CNP{C^\infty_{\vec{p_0}} (T^*Q)}$ is not only a subalgebra but also a
twosided ideal of $\CNP{C^\infty (T^*Q)}$. As we will see in
proposition \ref {WeylMainProp} it is {\em not} sufficient 
to consider only the functions with compact support in $T^*Q$.
For $f \in \CNP{C^\infty_{\vec{p_0}} (T^*Q)}$ we define
\BEQ {omegap}
    \omega_{\vec{p_0}} (f) :=
    \int_{\mathbb R^n} f(\vec q, \vec{p_0}) \, d^n q
\EEQ
which is obviously well-defined for any fixed momentum $\vec{p_0}$. 
Moreover it is a positive linear functional
of $\CNP{C^\infty_{\vec{p_0}} (T^*Q)}$:
\begin {proposition} \label {WeylPosProp}
Let $f, g \in \CNP{C^\infty_{\vec{p_0}} (T^*Q)}$ and define
$\Delta := \frac{\partial^2}{\partial q^k \partial p_k}$
and $S:= \exp(-\frac{i\lambda}{2} \Delta)$ then
\BEQ {omegaSS}
    \omega_{\vec{p_0}} (f*g) = \int_{\mathbb R^n}
    \left.\left(\left(S^{-1}f\right)\left(Sg\right)\right)
    \right|_{\vec{p} = \vec{p_0}} d^n q
\EEQ
and $\omega_{\vec{p_0}}$ is positive, i.e.
$\omega_{\vec{p_0}} (\cc f * f) \ge 0$.
\end {proposition}
\begin {proof}
A straight forward induction on $r$ shows that
\[
    r!\int_{\mathbb R^n} \left. \Lambda^r (f, g)
    \right|_{\vec{p} = \vec{p_0}} d^n q
    =
    \int_{\mathbb R^n} \left. \left(\sum_{s=0}^r {r \choose s}
    \Delta^s f (-\Delta)^{r-s} g \right)
    \right|_{\vec{p} = \vec{p_0}} d^n q
\]
then equation (\ref {omegaSS}) follows easily and this clearly implies
the positivity.
\end {proof}

Remark: A similar but non-formal positivity result of such integrals
was obtained in \cite [pp. 205-206 eqn. 5.11] {DFR95} in the context of
complex $C^*$-algebras.

Now we shall specialize to $\vec{p_0} = 0$ for simplicity.
We define for
$k = 1, \ldots, n$ the functions $P_k (\vec{q}, \vec{p}) := p_k$
and denote by $\mathcal J_0$ the left ideal of
$\CNP{C^\infty_Q (T^*Q)} := \CNP{C^\infty_{\vec{p_0} = 0} (T^*Q)}$
generated by the functions $P_k$ 
\BEQ {JDef}
    \mathcal J_0 := \left\{ f \in \CNP{C^\infty (T^*Q)} \; \left| \;
                    f = \sum\nolimits_k g_k * P_k, \;
                    g_k \in \CNP{C^\infty_Q (T^*Q)} \right\} \right.
\EEQ
which is indeed a left ideal of
$\CNP{C^\infty_Q (T^*Q)}$ since $\CNP{C^\infty_Q (T^*Q)}$ is a twosided
ideal of the whole algebra $\CNP{C^\infty (T^*Q)}$. By direct
calculation we get:
\begin {lemma} \label {JLem}
Let $f_1, \ldots, f_n \in \CNP{C^\infty_Q (T^*Q)}$ then we have
\BEQ {omeganull}
    \omega_0 \left( \cc{ (f_k * P_k)} * (f_l * P_l) \right) = 0
\EEQ
and $\mathcal J_0$ is contained in the Gel'fand ideal of $\omega_0$.
\end {lemma}
The next proposition will show that $\mathcal J_0$ is already equal to
the Gel'fand ideal of $\omega_0$ and the quotient
$\mathfrak H_0 := \CNP{C^\infty_Q (T^*Q)} / \mathcal J_0$ 
is canonically isomorphic to $\CNP{C^\infty_0 (Q)}$.
Let $r_0 : \CNP{C^\infty (T^*Q)} \to \CNP {C^\infty (Q)}$ be the
restriction to momentum $\vec{p} = 0$, i.e.
$r_0 (f) (\vec{q}) := f(\vec{q}, 0)$ and let $\pi_Q: T^*Q \to Q$ be the
projection to the configuration space. Then we clearly have
$r_0 (\CNP{C^\infty_Q (T^*Q)}) = \CNP{C^\infty_0 (Q)}$ and we define
$i_0 := \pi_Q^* \circ r_0$, i.e.
$i_0 (f)(\vec{q},\vec{p}) := f(\vec{q}, 0)$.
Moreover we need the operator
$I: \CNP{C^\infty_Q (T^*Q)} \to \CNP{C^\infty_Q (T^*Q)}$ defined for
$f \in \CNP{C^\infty_Q (T^*Q)}$ by
\BEQ {IDef}
    I(f) (\vec{q}, \vec{p}) := \int_0^1
    \frac{f(\vec{q}, t\vec{p}) - f(\vec{q}, 0)}{t} dt
\EEQ
and for $k = 1, \ldots, n$ we define
$T^k : \CNP{C^\infty_Q (T^*Q)} \to \CNP{C^\infty_Q (T^*Q)}$ by
\BEQ {TkDef}
    T^k (f) := \frac{\partial}{\partial p_k}
               \left( I \circ
               \frac{1}{1 + \frac{i\lambda}{2} \Delta \circ I} (f)
               \right) .
\EEQ
Note that $I(f), T^k (f) \in \CNP{C^\infty_Q (T^*Q)}$
if $f \in \CNP{C^\infty_Q (T^*Q)}$.
\begin {proposition} \label {WeylMainProp}
With the above notation we have for $f \in \CNP{C^\infty_Q (T^*Q)}$:
\begin {enumerate}
\item $f = i_0 \circ S (f) + T^k (f) * P_k$
\item This decomposition is unique, i.e.
      \BEQ {Zerlegung}
          \CNP{C^\infty_Q (T^*Q)} \cong
          i_0 \circ S \left(\CNP{C^\infty_Q (T^*Q)}\right)
          \oplus \mathcal J_0
      \EEQ
\item The Gel'fand ideal of $\omega_0$ is given by $\mathcal J_0$.
\end {enumerate}
\end {proposition}
\begin {proof}
First we use Hadamard's trick to obtain
$f = i_0(f) + P_k \partial_{p_k} I(f)$ which could be rewritten as
$f = i_0 (f) + (\partial_{p_k} I(f)) * P_k -
\frac{i\lambda}{2} \Delta \circ I(f)$ using the Weyl-Moyal product.
Iterating this equation leads to
$f = i_0 \circ \frac{1}{1 + \frac{i\lambda}{2} \Delta \circ I} (f)
+ T^k(f) * P_k$. An easy induction shows that
$r! i_0 \circ (\Delta \circ I)^r = i_0 \circ \Delta^r$ which
proves the first part. The second part follows from lemma \ref {JLem},
and the last statement is a consequence of part one and two.
\end {proof}

\begin {corollary}
The quotient space 
$\mathfrak H_0 := \CNP{C^\infty_Q(T^*Q)}/\mathcal J_0$ 
with the Hermitian product $\SP \cdot \cdot$ induced by 
$\omega_0$ is canonically isometric
to $\CNP{C^\infty_0 (Q)}$ with the Hermitian product
\BEQ {WeylHermProdDef}
    \SP \psi \phi := \int_{\mathbb R^n}
                     \cc {\psi (\vec{q})} \phi (\vec{q}) \, d^n q
    \qquad \qquad
    \psi, \phi \in \CNP{C^\infty_0 (Q)}
\EEQ
and the canonical isometry $r$ is induced by $r_0$:
For $\psi_f \in \mathfrak H_0$ we set $r(\psi_f) := r_0 \circ S(f)$ and
the inverse of $r$ is simply the pull-back ${\pi_Q^* \bmod \mathcal J_0}$.
\end {corollary}

Now we consider the GNS representation $\pi_0$ induced by $\omega_0$ on
$\mathfrak H_0$. First we notice that according to corollary
\ref {GNSIdealRepCor} not only
$\CNP{C^\infty_Q (T^*Q)}$ could be represented but the whole algebra
$\CNP{C^\infty (T^*Q)}$ since $\CNP{C^\infty_Q (T^*Q)}$ is a twosided
ideal and $\mathcal J_0$ is a left ideal of $\CNP{C^\infty (T^*Q)}$.
Let $\varrho$ be the corresponding representation on
$\CNP{C^\infty_0 (Q)}$, i.e.
$\varrho (f) := r \circ \pi_0 (f) \circ \pi^*_Q$.
\begin {theorem} [Formal Schr\"odinger Quantization]
Let $Q = \mathbb R^n$ and $\psi \in \CNP{C^\infty_0 (Q)}$ be a
`formal wave function' on the configuration space $Q$ and
$f \in \CNP{C^\infty (T^*Q)}$.
\begin {enumerate}
\item \BEQ {SchrRep}
          \varrho(f) \psi = \sum_{r=0}^\infty \frac{1}{r!}
          \left(\frac{\lambda}{i}\right)^r
          \left.\frac{\partial^r (Sf)}
          {\partial p_{i_1} \cdots \partial p_{i_r}}
          \right|_{\vec{p} = 0}
          \frac{\partial^r \psi}
          {\partial q^{i_1} \cdots \partial q^{i_r}}
      \EEQ
\item For polynomials in $q^1, \ldots, p_n$ the representation $\varrho$
      is the canonical quantization rule, i.e.
      \BEQ {CanQRule}
          \varrho (q^k) = q^k \qquad
          \varrho (p_k) = \frac{\lambda}{i}
          \frac{\partial}{\partial q^k}
      \EEQ
      and the polynomials are mapped to the Weyl ordered polynomials of
      the corresponding operators (symmetrization rule).
\item The $\CNPC$-pre-Hilbert space $\CNP{C^\infty_0 (Q)}$ with the
      Hermitian product (\ref {WeylHermProdDef}) is already Cauchy
      complete and hence a $\CNPC$-Hilbert space.
\item $\mathfrak H_0$ does not admit a (countable or uncountable)
      Hilbert base (for $n \ge 1$).
\end {enumerate}
\end {theorem}
\begin {proof}
The first statement is a straight-forward computation. For the second
one notes first that for each nonnegative integer $k$ and $2n$ formal
parameters $\alpha_1,\ldots,\alpha_n,\beta^1,\ldots,\beta^n$ the function
$(\alpha_rq^r+\beta^rp_r)^k$ is assigned the operator
$(\alpha_r \varrho(q^r)+\beta^r\varrho(p_r))^k$ by the Weyl
symmetrization rule. Hence the formal exponential function
$e_{\alpha,\beta}(\vec{q},\vec{p}):=\exp(\alpha_rq^r+\beta^rp_r)$ is assigned
the operator $\exp(\alpha_r\varrho(q^r)+\beta^r\varrho(p_r))$. On the other
hand, the right hand side of (\ref{SchrRep}) is easily seen to be equal to
the standard ordering prescription of the function $Sf$, i.e. where after
applying the rule (\ref{CanQRule}) all the 
derivatives are put on the right hand side first. Clearly, the function
$Se_{\alpha,\beta}=\exp(i\lambda\alpha\beta/2)e_{\alpha,\beta}$ is assigned
the operator
$\exp(i\lambda\alpha\beta/2)\exp(\alpha_i\varrho(q^i))
                              \exp(\beta^i\varrho(p_i))$ 
by standard ordering.
Now the Baker-Campbell-Hausdorff formula for the $2n+1$-dimensional
Heisenberg Lie algebra easily implies that the Weyl ordered operator for
$e_{\alpha,\beta}$ and the standard-ordered operator for $Se_{\alpha,\beta}$
coincide which proves this statement.
The third
and fourth statement are proved in a more general context in the next
theorem.
\end {proof}

\begin {theorem}
Let $M$ be an orientable manifold with volume form $\Omega$. Then
$\mathfrak H := \CNP{C^\infty_0 (M)}$ together with the integral
\[
    \SP f g := \int_M \cc f g \, \Omega \qquad f, g \in \mathfrak H
\]
as Hermitian product is a $\CNPC$-Hilbert space.
Moreover $\mathfrak H$ does not admit a
(countable or uncountable) Hilbert base if $\dim M \ge 1$.
\end {theorem}
\begin {proof}
 It is easy to see that every Cauchy sequence in
 $\mathfrak H$ with respect to the norm induced by the scalar product
 is a Cauchy sequence in the metric sense of Proposition \ref{CNPCompProp},
 and vice versa, which proves completeness.

 Suppose there would exist a Hilbert base
 $\{\psi^{(\alpha)}\}_{\alpha\in I}$,
 $I$ some index set, for $\mathfrak H$. Because of its orthonormality
 every element is of the form $\psi^{(\alpha)}=\psi^{(\alpha)}_0+
 \lambda^{a_{\alpha}}\psi^{(\alpha)}_1$ where the $\psi^{(\alpha)}_0$ is
 an
 orthonormal family in $C^\infty_0(M)$, $a_\alpha$ is a positive rational
 number, and $\psi^{(\alpha)}_1\in\CNP{C^\infty_0(M)}$ with vanishing
 coefficients of negative $\lambda$-powers. By the usual complex
 $L^2$-theory for manifolds of dimension greater or equal than 1 there
 are countably many different $\psi^{(\alpha)}_0$, and there is a function
 $f\in C^\infty_0(M)$ which is not a finite linear combination of the
 $\psi^{(\alpha)}_0$. If $f$ were approximated by finite sums of the
 $\psi^{(\alpha)}$ there would be a positive integer $N$ and $\alpha_1,
 \ldots,\alpha_N\in I$ such that
 \[
    \lambda^2 > \norm{f-\sum_{i=1}^N\SP{\psi^{(\alpha_i)}}{f}}^2
                 = \SP{f}{f}-\sum_{i=1}^N\SP{f}{\psi^{(\alpha_i)}}
                                         \SP{\psi^{(\alpha_i)}}{f}
 \]
 In particular it would follow that
 $\SP{f}{f}=\sum_{i=1}^N\SP{f}{\psi^{(\alpha_i)}_0}
                                         \SP{\psi^{(\alpha_i)}_0}{f}$
 implying by Parseval's equality for complex Hilbert spaces that $f$ is
 a finite linear combination of the $\psi^{(\alpha_i)}$ which is a
 contradiction.
\end {proof}

At last we want to describe again the way back to convergence if we
substitute the formal parameter $\lambda$ by $\hbar \in \mathbb R^+$ and
again we want to ask for convergence in the representation.
Let $\psi = \sum_{r=0}^\infty \lambda^{q_r} \psi_{q_r}
\in \CNP{C^\infty_0 (Q)}$ be a formal wave function and define the 
truncated series
$\psi_N := \sum_{r=0}^N \lambda^{q_r} \psi_{q_r}$. Then we ask for
convergence in the sense of the locally convex topology of the 
smooth functions with compact support $\mathcal D (Q) := C^\infty_0 (Q)$ 
of the sequence $\widetilde{\psi}_N$ if we substitute $\lambda$ by 
$\hbar$ and define
\begin {eqnarray*}
    H(\hbar)
    & := & \left\{
           \psi \in \CNP{C^\infty_0(\mathbb R^n)} \; \left| \;
           \psi_N|_{\lambda = \hbar} 
           \stackrel{\mathcal D}{\longrightarrow}
           \Psi \in \mathcal D (Q) 
           \mbox { as } N \to \infty \right\}
           \right. \label {WeylHhbarDef} \\
    N(\hbar)
    & := & \left\{ \psi \in H(\hbar) \; \left| \;
           \psi_N|_{\lambda = \hbar} 
           \stackrel{\mathcal D}{\longrightarrow} 0
           \mbox { as } N \to \infty
           \right\} \right. \label {WeylNhbarDef} .
\end {eqnarray*}
Now the same procedure as in the Wick case can be done and we get 
using the well-known completeness 
(see e.~g. \cite [Theorem 6.5.(g)]{Rud91}) of $\mathcal D(Q)$:
\begin {theorem}
With the notations of above we have:
\begin {enumerate}
\item $H(\hbar)$ is a $\mathbb C$-vector space and $N(\hbar)$ is a
      $\mathbb C$-subvector space of $H(\hbar)$.
\item If $\psi, \phi \in H(\hbar)$ then $\SP \psi \phi$ converges for
      $\lambda = \hbar$ and
      \[
          \SP \psi \phi \big|_{\lambda = \hbar}
          = \int_{\mathbb R^n} \cc {\Psi(\vec q)} \Phi (\vec q) \, d^n q
      \]
      and defines a positive semidefinite sesquilinear form for
      $H(\hbar)$.
\item The quotient $\mathcal H(\hbar) := H(\hbar)/N(\hbar)$
      is canonically isometric to the $\mathbb C$-pre-Hilbert space
      $\mathcal D (Q)$ where $\psi \bmod N(\hbar) \mapsto \Psi$ 
      is the isomorphism.
\end {enumerate}
\end {theorem}
Then the usual completion of $\mathcal D (Q)$ 
with respect to the $L^2$-norm leads to the Hilbert space 
$L^2 (\mathbb R^n)$. Again the GNS representation $\varrho$ induces a 
representation of at least `many' elements 
$f \in \CNP{C^\infty (\mathbb R^{2n})}$
on $\mathcal H(\hbar) \cong  \mathcal D (Q)$. For
$f \in \CNP{C^\infty (\mathbb R^{2n})}$ such that
$\varrho(f) N(\hbar) \subseteq N(\hbar)$ we define again $D_f (\hbar)$
and $\mathcal D_f (\hbar)$ as in the Wick case and obtain a quantization
map $\mathcal Q$ defined by
\BEQ {WeylQuant}
    \mathcal Q(f) (\psi \bmod N(\hbar))
    := \varrho (f) \psi \bmod N(\hbar).
\EEQ
for $\psi \in D_f (\hbar)$. Then $\mathcal Q(f)$ is a linear operator
defined on some domain $\mathcal D_f (\hbar)
\subseteq \mathcal H(\hbar) \cong \mathcal D(\mathbb R^n)$ and again 
there are `many' elements $f$ with dense domain $\mathcal D_f (\hbar)$,
for example those functions which are polynomial in the momentum 
variables which have $\mathcal D_f (\hbar) = \mathcal D(Q)$ 
and in particular the polynomials in $\vec q$ and $\vec p$.
In this case the polynomials are represented via $\mathcal Q$ 
by differential operators obtained form (\ref {SchrRep}) after 
substituting $\lambda$ by $\hbar$.

Note that an analogous result holds if one substitutes `$\mathcal D$' by 
`$\mathcal S$' denoting Schwartz's test function space 
with its locally convex topology resp. by `$L^2$' denoting 
the usual square integrable functions. In the case of e.~g. polynomials 
in $\vec q$ and $\vec p$ the above quantization map is the 
{\em Weyl transformation} which is discussed e.~g. in 
\cite {Kam86,Mai86} where among other things the image of the Weyl 
transformation applied to certain tempered distributions on 
$\mathbb R^{2n}$ is computed.

\section{Open Problems}
\label {OutSec}

In this section we list some open problems arising with our approach:
\begin{enumerate}
\item It would be very interesting to see whether the approach 
      of the WKB approximation via Lagrangean submanifolds
      (see e. g. \cite {BW95}) contained in the energy
      surface of a fixed Hamiltonian function on the 
      symplectic manifold is related to some variant of GNS 
      construction with respect to a suitable positive linear 
      functional whose support is contained in
      that Lagrangean submanifold. In the particular case 
      of an arbitrary projectable Lagrangean submanifold of 
      $T^*\mathbb R^n$ we have been able to derive a 
      suitable GNS construction for the usual WKB expansion in 
      \cite {BW96b}.       
      
      A further problem arising in this context is the
      fact that the discrete energy eigenvalues in quantum mechanics are
      dependent on $\hbar$, and it is not obvious how this can 
      be encoded in a formal way in the energy surface.
\item We have seen in section \ref {WickSec} that the spectrum 
      of the harmonic oscillator can be computed in a purely 
      formal manner before the convergence
      scheme is performed. It may be interesting to develop a kind of
      `formal spectral theory' in order to formally compute e.g.
      discrete spectra
      depending on $\lambda$ and then deal with the convergence.
\item Related to this question one may ask more generally to what
      extend there is some reasonable functional analysis in these formal
      Hilbert spaces: i.e. which (possibly weaker) topologies are more
      suited for the definition and calculation of spectra and convergence
      properties.                               
      The usual literature on $p$-adic functional analysis
      (see e.g. \cite{NBB71,Kal47}) and $p$-adic formulations 
      of quantum physics (see e.g. \cite {FO87,FW87,AV91,AK96})
      does unfortunately not seem to 
      deal with Hilbert spaces in the sense described in our 
      approach, but exclusively with fields carrying absolute values 
      such as the field $\mathbb Q_p$ of $p$-adic numbers (which is 
      non-orderable, see appendix A), and seems to avoid
      ordered fields and hence the notion of positivity we need 
      (cf. also the remark in \cite [p. 5517] {AK96}).
\item Finally, it may be interesting to see whether the prequantum line
      bundles of geometric quantization over, say, a compact prequantizable
      K\"ahler manifold (see e.g. \cite{Raw77}, \cite{CGRI}, \cite{CGRII})
      are related to this construction (perhaps after employing a
      convergence scheme) and if yes, whether they can be constructed that
      way.

\end{enumerate}

\appendix

\section {Some analytic properties of ordered fields}
\label {AppA}

In this appendix we shall examine some standard definitions of
topology and calculus and transfer them from the case of real numbers
$\mathbb R$ to the more general case of an arbitrary ordered field
$\field R$. Again all the statements can be proved analogously using
only the ordering axioms.

First of all we define $\epsilon$-balls around any point in $\field R$ 
by means of the ordering relation and notice that the ordering induces
a {\em topology} and a {\em uniform structure} for $\field R$ similar
to the topology and the uniform structure of a metric space
\cite [Chap. 6] {Kel55}.
\begin {definition}
Let $\field R$ be an ordered field. For any
$x \in \field R$ and any $0 < \epsilon \in \field R$ we define the
$\epsilon$-ball around $x$ by
\BEQ {eBall}
    B_\epsilon(x) := \left\{ y \in \field R \; \big| \;
                     |x-y| < \epsilon \right\} .
\EEQ
The set of all $\epsilon$-balls is denoted by
$\mathcal B := \left\{ B_\epsilon (x) \subseteq \field R \; | \;
0 < \epsilon \in \field R, x \in \field R \right\}$.
Furthermore we define for $0 < \epsilon \in \field R$
\BEQ {eStreifen}
    U_\epsilon := \left\{ (x,y) \in \field R \times \field R \; \big| \;
                  |x-y| < \epsilon \right\}
                  \subseteq \field R \times \field R
\EEQ
and
$\mathcal U' := \left\{ U_\epsilon \subseteq \field R \times \field R
  \; | \; 0 < \epsilon \in \field R \right\}$.
\end {definition}
\begin {proposition} \label {TopUniProp}
The set of $\epsilon$-balls $\mathcal B$ is a base of a topology
$\mathcal T$ for the field $\field R$ such that $\field R$ becomes a
normal space. The set $\mathcal U'$ is a base for a uniform
structure $\mathcal U$ for $\field R$ which induces the same topology
$\mathcal T$ than $\mathcal B$. Addition and multiplication are
uniformly continuous maps from $\field R \times \field R \to \field R$.
\end {proposition}
We shall call this topology and the uniform structure
induced by the ordering relation the standard topology and the standard
uniform structure of the ordered field $\field R$. Using this topology
we can define continuous functions $f:\field R \to \field R$ as usual.

Remark: We notice that in the
topology induced by the ordering relation the intervals may neither be
connected nor compact in general.
As an example we consider the field of formal Laurent series $\LR$ over
the real numbers and the open balls
$B_{n\lambda} (r)$ with $r \in \mathbb R$ and $n \in \mathbb N$.
Then for $a<b$ the closed interval
$[a,b] := \{x | a \le x \le b\} \subset \LR$ in the Laurent field is
contained in
\[
    [a,b] \subset \bigcup_{{n\in \mathbb N \atop r \in [a,b] \subset
    \mathbb R}} B_{n\lambda} (r)
\]
but $B_{n\lambda} (r) \cap B_{n'\lambda} (r') = \emptyset$ if
$r \ne r'$. Hence $[a,b] \subset \LR$ is neither compact nor
connected.

Another important concept is the supremum and infimum of a bounded
subset in an ordered field. We define them as usual:
A subset $U \subseteq \field R$ is called {\em bounded} iff there is
a $C \in \field R$ such that $|x| \le C$ for all $x \in U$.
Then $C$ is called a bound for $U$. Analogously we define upper
and lower bounds. $C$ is called the {\em supremum (infimum)}
of $U$ iff $C$ is the smallest upper bound (largest lower bound),
i.e. for all $0<\epsilon \in \field R$ there exists
$x \in U$ such that $x > C - \epsilon$ ($x < C + \epsilon$).
In an archimedian ordered field the existence of a supremum or
infimum of a bounded set can be proved iff the field is Cauchy complete.
On the other hand, in a non-archimedian ordered field
neither the supremum nor the infimum of a bounded set exist in general
\cite [p. 245] {vdWI}.

\begin {definition} [Convergence, Cauchy sequences, completeness]
Let $\field R$ be an ordered field and $(a_n)_{n\in \mathbb N}$ a
sequence in $\field R$.
Then $(a_n)$ is called convergent to $a \in \field R$ iff
\[
    \forall \, 0<\epsilon \in \field R \quad
    \exists N \in \mathbb N
    \quad \mbox{\it such that} \quad \forall n>N : |a_n - a| < \epsilon .
\]
A sequence $(a_n)$ is called a Cauchy sequence iff
\[
    \forall \, 0<\epsilon \in \field R \quad
    \exists N \in \mathbb N
    \quad \mbox{\it such that} \quad \forall m,n>N : 
    |a_m - a_n| < \epsilon .
\]
An ordered field is called Cauchy complete iff every Cauchy sequence
converges in $\field R$.
\end {definition}
It is a well-known result that for any ordered field $\field R$ there 
exists a unique ordered field $\widehat{\field R}$ such that
$\field R$ is a dense subfield of $\widehat{\field R}$, the
orderings are compatible and $\widehat{\field R}$ is Cauchy
complete \cite [p. 238] {vdWI}.

For certain approximation statements we shall need to know whether the
topology or
the uniform structure is first countable or not
(a topology for a set
is called first countable iff any point has a countable base of
neighbourhoods \cite [p. 50] {Kel55}).
In our case the existence of such a countable
base of neigbourhoods is equivalent to the existence of
non-trivial zero-sequences:
\begin {proposition} \label {FirstCountProp}
Let $\field R$ be an ordered field. Then the following 
properties are equivalent:
\begin {enumerate}
\item There is a sequence $(\epsilon_n)_{n \in \mathbb N}$ such that
      $\epsilon_n > 0$ for all $n \in \mathbb N$ and 
      $\epsilon_n \to 0$.
\item The standard uniform structure $\mathcal U$ has a countable base.
\item The standard topology $\mathcal T$ is first countable.
\item The standard uniform structure $\mathcal U$ and the standard
      topology $\mathcal T$ can be induced by a metric.
\end {enumerate}
\end {proposition}
\begin {proof} The implications
$i.) \Rightarrow ii.) \Rightarrow iii.)$ are obvious and
$iii.) \Rightarrow i.)$ is proved by constructing the zero sequence
using a given countable neighbourhood base. The fourth part is
equivalent to the second one for general reasons since
$\mathcal T$ is Hausdorff according to proposition \ref {TopUniProp}
(cf. e.g. \cite[p. 186]{Kel55}).
\end {proof}

Although the supremum of a bounded subset does not exist in general
we can easily prove the existence of a sequence in a bounded subset with
supremum which converges to the supremum if the field is first countable:
\begin {lemma} \label {ConvSupLem}
Let $\field R$ be an ordered field such that the standard 
topology is first countable. If for 
$A \subset \field R$ the supremum $\sup A$ exists then there is a 
sequence $(a_n)$ with elements in $A$ such that $a_n \to \sup A$.
\end {lemma}

In the theory of fields 
another possibility to define a metric topology is an
{\em absolute value}:
\begin {definition} [{Absolute value \cite [p. 558] {JacI}}]
Let $\field R$ be a field. An absolute value $\varphi$ is a map
$\varphi : \field R \to \mathbb R^+\cup \{0\}$ such that for all
$a,b \in \field R$
\begin {enumerate}
\item $\varphi(a) = 0 \Longleftrightarrow a = 0$
\item $\varphi(ab) = \varphi(a)\varphi(b)$
\item $\varphi(a+b) \le \varphi(a) + \varphi(b)$ .
\end {enumerate}
An absolute value $\varphi$ is called non-archimedian iff
$\varphi(a+b) \le \max (\varphi(a), \varphi(b))$
and archimedian if this is not the case. 
An absolute value is called trivial iff 
$\varphi (0) = 0$ and $\varphi (a) = 1$ for all $0 \ne a \in \field R$.
\end {definition}
With help of an absolute value one can define a metric on $\field R$
for $a,b \in \field R$ by
\BEQ {ValMetricDef}
    d_\varphi (a,b) := \varphi (a-b) .
\EEQ
Recall that if the absolute value is
non-archimedian then the metric $d_\varphi$ is an {\em ultrametric},
i.e. $d_\varphi (a,b) \le \max(d_\varphi(a,c), d_\varphi(c,b))$
\cite [p. 6] {NBB71}.
If $\varphi: \field R \to \mathbb R$ is an absolute value we denote
by $\mathcal T_\varphi$ and $\mathcal U_\varphi$ the 
topology and the uniform structure induced by the corresponding metric
$d_\varphi$. The open metric balls around $x \in \field R$ with radius
$0 < \epsilon \in \mathbb R$ are denoted by $B_\epsilon^\varphi (x)$ and
analogously we define
$U_\epsilon^\varphi := \{ (x,y) \in \field R \times \field R
\; | \; d_\varphi (x,y) < \epsilon \}$ for $0 < \epsilon \in \mathbb R$.
Then the metric balls $B_\epsilon^\varphi (x)$ and the
$U_\epsilon^\varphi$ form a base for the topology $\mathcal T_\varphi$
and the uniform structure $\mathcal U_\varphi$.
As in the case of an ordered field the field $\field R$ with 
absolute value $\varphi$ can be completed with respect to $d_\varphi$
in the usual sense of metric completion.

Remark: It is well-known that the field $\mathbb Q_p$ of $p$-adic 
numbers (which is the metric completion of the field $\mathbb Q$ 
of rational numbers by means of the $p$-adic absolute value 
\cite [p.558] {JacII}) is not orderable: note that the series
$z:=1/\sqrt{1-4p}:=\sum_{i=0}^\infty{2k \choose k}p^k$ converges in 
$\mathbb Q_p$ and $-1=(4p-1)z^2$ hence $-1$ is a sum of squares
in $\mathbb Q_p$.

In the case of an ordered field $\field R$ we ask now for an absolute
value $\varphi$ such that $\varphi: \field R \to \mathbb R$ is not only
continuous with respect to $\mathcal T$
but such that $\mathcal T_\varphi = \mathcal T$ and
$\mathcal U_\varphi = \mathcal U$. If this is the case we call the
absolute value {\em compatible} with the ordering of $\field R$.
\begin {lemma} \label {TopUnifLem}
Let $\field R$ be an ordered field and
$\varphi : \field R \to \mathbb R$ an absolute value. Then
$\mathcal T = \mathcal T_\varphi \iff \mathcal U = \mathcal U_\varphi$.
\end {lemma}
\begin {proof}
The inverse implication being trivial assume that $\mathcal T =\mathcal
T_{\varphi}$. Then for all $0<\epsilon\in\mathbb R$ there exists
$0<\delta\in\field R$ such that $B_\delta(0)\subset
B_\epsilon^\varphi(0)$. This implies $U_\delta\subset U_\epsilon^\varphi$
and vice versa.
\end {proof}

If there is such a compatible absolute value then the
topology $\mathcal T = \mathcal T_\varphi$ of $\field R$ is first
countable since it is induced by a metric and hence we can apply
proposition \ref {FirstCountProp}. Furthermore $\varphi$ is continuous
which implies that $\varphi$ is non-trivial.
This implies that there are elements $\epsilon \in \field R$ such that
$0 < \varphi(\epsilon) < 1$ and hence $\varphi (\epsilon)^n \to 0$.
Thus $\epsilon^n \to 0$ in the field $\field R$ since the
topologies $\mathcal T_\varphi$ and $\mathcal T$ are the same.
Hence we have proved the following lemma:
\begin {lemma}
Let $\field R$ be an ordered field and
$\varphi : \field R \to \mathbb R$ a compatible absolute value. Then
there exists $0 < \epsilon \in \field R$ such that
$\epsilon^n \to 0$.
\end {lemma}

At last we shall consider ordered fields and their real closure.
A field $\field R$ is called {\em real closed} iff it is ordered
and any positive element has a
square root in $\field R$ and every polynomial of odd degree with
coefficients in $\field R$ has a root in $\field R$ \cite [p. 308] {JacI}. 
Let $\field R$ be an ordered field. An extension field
$\widehat{\field R}$ of
$\field R$ is called real closure of $\field R$ iff $\widehat{\field R}$
is real closed and algebraic over $\field R$ and the (unique) order
of $\widehat{\field R}$ is an extension of the order in $\field R$
\cite [p. 655] {JacII}. Such a real closure exists for every orderd
field and is unique up to isomorphisms \cite[p. 656]{JacII}.
Furthermore for a real closed field one can prove the
`fundamental theorem of algebra': If $\field R$ is a real closed field
then the quadratic field extension
$\field C = \field R(i)$ with $i^2 := -1$ is algebraically
closed \cite [p. 309] {JacI}.

For the quadratic field extension
$\field C = \field R(i)$ of a real closed field
$\field R$ we define
\BEQ {zBetrag}
    |z| := \sqrt{z\cc z} = \sqrt{a^2 + b^2}.
\EEQ
for $z = a + ib \in \field C$ where $a,b \in \field R$ using the square
root. Then $|z| \in \field R$ and for the real elements $a \in \field R$
this definition coincides with the previous definition since
$\sqrt{a^2} = |a|$. Furthermore we have $|z| = 0 \iff z = 0$,
$|zw| = |z||w|$ and $|z+w| \le |z| + |w|$ for $z,w \in \field C$
and hence $|\cdot|$ induces a topology on $\field C$ which is compatible
with the standard topology of the ordered field $\field R$.
Again we will call this topology for $\field C$ the standard topology.
Then the field $\field R$ is Cauchy complete with respect to the
standard topology iff its quadratic extension $\field C$ is Cauchy
complete with respect to the standard topology.

\section {Hilbert space formulation based on ordered fields
          and the generalized GNS construction}
\label {AppB}

In this appendix we shall give a generalization of the usual definition
of a Hilbert space over $\mathbb C$ which leads to the definition of
a Hilbert space over $\field C = \field R(i)$ where
$\field R$ is an ordered field which is Cauchy complete {\em and} real
closed and formulate the generalized GNS construction based 
thereon.

Let $\HC$ be a pre-Hilbert space
over $\field C = \field R(i)$ where $\field R$ is
a real closed field. In this case we can define a
{\em $\field R$-norm} for vectors in $\HC$
that takes values in $\field R$ by
\BEQ {KNormDef}
    \norm{\phi} := \sqrt{\SP \phi \phi}
\EEQ
for $\phi \in \HC$. Then we clearly have for
$\phi, \psi \in \HC$ and $a \in \field C$:
\BEQ {KNormProps}
    \begin {array} {c}
    \norm{\phi} \ge 0 \quad
    \mbox{\rm and } \quad \norm{\phi} = 0 \iff \phi = 0 \\
    \norm {a\phi} = |a| \norm{\phi} \\
    \norm {\phi + \psi} \le \norm{\phi} + \norm {\psi} \\
    \norm {\phi + \psi}^2 + \norm {\phi - \psi}^2
    = 2\norm{\phi}^2 + 2 \norm{\psi}^2
    \end {array}
\EEQ
With this $\field R$-norm the pre-Hilbert space
$\HC$ becomes a topological vector space where we
define a base for the topology by the $\epsilon$-balls with respect to
$\norm\cdot$. Note that $\HC$ is {\em not} a normed
vector space in the usual sense since $\norm\cdot$ takes values
in $\field R$ and not in $\mathbb R$ (in contrast, for example, 
to Kalisch's definition of $p$-adic Hilbert 
spaces \cite [p. 181] {Kal47}). Note also that in order to define the 
Hilbert norm (\ref {KNormDef}) we need the existence of square 
roots of every element of $\field C$ which motivates the use of 
algebraically closed fields.

An easy consequence is the
following lemma:
\begin {lemma} \label {SPStetigLem}
Let $\HC$ be a pre-Hilbert space over
$\field C = \field R(i)$ where $\field R$ is a
real closed field. Then the Hermitian product
$\SP \cdot \cdot : \HC \times \HC \to \field C$ is continuous with
respect to $\norm \cdot$.
\end {lemma}
By means of $\norm\cdot$ we can define Cauchy sequences in
$\HC$ and this leads to the definition of a Hilbert space:
\begin {definition} [Hilbert space]
Let $\field C = \field R(i)$ be the quadratic field
extension of a real closed field $\field R$ and $\HC$ a pre-Hilbert
space over $\field C$ with $\field R$-norm $\norm\cdot$. Then
$\HC$ is called a Hilbert space over $\field C$ iff every Cauchy sequence
in $\HC$ converges in $\HC$ with respect to $\norm\cdot$.
\end {definition}
\begin {lemma}
Let $\{0\} \ne \HC$ be a Hilbert space over
$\field C = \field R(i)$
where $\field R$ is a real closed field. Then $\field C$
and $\field R$ are Cauchy complete.
\end {lemma}
Due to this lemma we should always assume that $\field C$ is not
only algebraically closed but also Cauchy complete if we want to
consider Hilbert spaces over $\field C$. In a next step we want to
construct the completion of a pre-Hilbert space to a Hilbert space. 
\begin {proposition}
\label {CompleteHilbProp}
Let $\HC$ be a pre-Hilbert space over 
$\field C = \field R(i)$ where $\field R$ is real closed
and Cauchy complete. Then up to unitary equivalence there is one
Hilbert space $\widehat{\HC}$ such that there is an
isometry $i: \HC \to \widehat{\HC}$ and $i(\HC)$ is dense in
$\widehat{\HC}$.
\end {proposition}
Let $\HC$ and $\KC$ be Hilbert spaces over $\field C$ and
$T:\HC \to \KC$ a linear map. Then clearly $T$ is continuous iff it is 
continuous at some point $\phi_0 \in \HC$ since $T$ is linear.
\begin {lemma}
Let $\HC$ and $\KC$ be Hilbert spaces over 
$\field C = \field R(i)$ where $\field R$ is a real closed
and Cauchy complete field. For a linear map $T: \HC \to \KC$ we have:
$T$ is continuous iff there exists $C \in \field R$ such that 
$\norm{T\phi} \le C \norm{\phi}$ for all $\phi \in \HC$ ($T$ is bounded).
\end {lemma}
Remark:
For continuous linear maps $T$ between Hilbert spaces over $\mathbb C$
one usually defines an operator norm by the supremum of
$\norm {T\phi}/\norm{\phi}$ where $\phi \ne 0$ ranges over the Hilbert
space. But in the case of a non-archimedian ordered field such 
a supremum does not exist in general though $T$ is bounded.

In the case that the standard topology of the
field $\field R$ is first countable we have the following important
property of $\field C$-Hilbert spaces:
\begin {lemma} \label {WDichtLem}
Let $\HC$ be a $\field C$-Hilbert space
where $\field C = \field R(i)$ and $\field R$ is real
closed, Cauchy complete and first countable. Let $W$ be a subspace
of $\HC$. If $W$ is dense in $\HC$ then for any
$\phi \in \HC$ there exists a sequence
$\phi_n \in W$ with $\phi_n \to \phi$.
\end {lemma}

An important concept in the usual theory of Hilbert spaces are the
Hilbert bases. This leads to the following generalizations:
Let $\HC$ be a $\field C$-Hilbert space
where $\field C = \field R(i)$ and $\field R$ is a real
closed and Cauchy complete field. Let $\mathbb I$ be an index set
and let $\cal F$ be the set of all finite subsets of $\mathbb I$.
Let $\{ e_k\}_{k \in \mathbb I}$ be a set of
vectors in $\HC$ 
then $\{ e_k \}_{k\in \mathbb I}$ is called an {\em orthonormal system}
iff $\SP {e_k} {e_{k'}} = \delta_{kk'}$ for all $k, k' \in \mathbb I$. 
In this case the set $\{ e_k \}_{k \in \mathbb I}$ is clearly 
linear independent in $\HC$. This leads to the definition 
of a Hilbert base:
\begin {definition} [Hilbert base]
Let $\HC$ be a Hilbert space and let
$\{ e_k \}_{k \in \mathbb I}$ be an orthonormal set
where $\field C = \field R(i)$ and $\field R$ is real
closed, Cauchy complete and first countable.
Then $\{e_k\}_{k\in \mathbb I}$ is called a Hilbert base for $\HC$ iff
$\field C$-span$\{e_k\}_{k\in \mathbb I}$ is dense in $\HC$ with respect
to the topology induced by $\norm\cdot$.
\end {definition}
The following theorem is proved in a completely analogous fashion to the
usual proofs in textbooks on functional analysis (see e.g.
\cite[p. 86]{Yos80})
\begin {theorem}
Let $\HC$ be a $\field C$-Hilbert space where
$\field C = \field R(i)$ and $\field R$ is a real closed
and Cauchy complete field. Let $\{ e_k \}_{k \in \mathbb I}$ be an
orthonormal set. Then we have {\em Bessel's inequality}
\BEQ {BesselIQ}
    \sum_{k\in F} \left| \SP {e_k} \phi \right|^2
    \; \le \; \SP \phi \phi
    \qquad \qquad
    \forall \phi \in \HC, \forall F \in {\cal F}
\EEQ
and for all $\phi \in \HC$
the following {\em `best approximation property'}
\BEQ {BestApp}
    \norm {\phi - \sum_{k\in F} \alpha_k e_k}
    \; \ge \;
    \norm {\phi - \sum_{k\in F} \SP {e_k} \phi e_k}
\EEQ
for all $\alpha_k \in \field C$ where the equality is only
satisfied iff $\alpha_k = \SP{e_k}\phi$.
Now let in addition $\field R$ be first countable, let
$\mathbb I = \mathbb N$, and let
$\{ e_k \}_{k \in \mathbb N}$ be a Hilbert base for $\HC$.
Then there exists a sequence
$\phi_n = \sum_{k=1}^{N_n} \alpha^n_k e_k
\in \field C$-span$\{e_k\}_{k \in \mathbb I}$
with $\phi_n \to \phi$ for any $\phi \in \HC$ and
\BEQ {ekTophi}
     \phi = \lim_{N\to \infty} \sum_{k=1}^N \SP{e_k}\phi e_k
     =: \sum_{k=1}^\infty \SP{e_k}\phi e_k
\EEQ
and we have {\em Parseval's equation}:
\BEQ {Parseval}
    \norm {\phi}^2 = \sum_{k=1}^\infty \left|\SP {e_k} \phi\right|^2
    \qquad \qquad
    \SP \phi \psi =
    \sum_{k=1}^\infty \SP \phi {e_k} \SP {e_k} \psi
\EEQ
Moreover if $\alpha_k \in \field C$ then $\sum_{k} \alpha_k e_k$
converges to a vector in $\HC$ iff $\sum_{k} |\alpha_k|^2$
converges in $\field R$.
An orthonormal system $\{e_k\}_{k\in \mathbb I}$ is a Hilbert base
for $\HC$ iff Parseval's equation holds for any $\phi \in \HC$.
\end {theorem}
As a (generic) example we will consider the $\ell^2$-space of a
field $\field C = \field R(i)$ where $\field R$ is real
closed, Cauchy complete and first countable:
\BEQ {LzweiDef}
    \ell^2 (\field C) :=
    \left\{ (a_k)_{k \in \mathbb N} \; \left| \;
    a_k \in \field C \mbox { such that }
    \sum_{k=1}^\infty |a_k|^2
    \mbox { converges in } \field R \right\} \right.
\EEQ
together with the $\ell^2$-product
\BEQ {LzweiProdDef}
          \SP a b :=
          \sum_{k=1}^\infty \cc a_k b_k .
\EEQ
\begin {proposition} \label {EllzweiProp}
Let $\field C = \field R(i)$ be the quadratic field
extension of $\field R$ where $\field R$ is real closed, Cauchy complete
and first countable. Then $\ell^2(\field C)$ is a
$\field C$-Hilbert space and the vectors
$\hat e_k := (0, \ldots, 0, 1, 0, \ldots)$ (where $1$
is the $k$-th entry) form a countable
Hilbert base for $\ell^2(\field C)$.
\end {proposition}
\begin {proposition} \label {UniEquiProp}
Let $\HC$ and $\KC$ be Hilbert spaces
where $\field C = \field R(i)$  and $\field R$ is real
closed, Cauchy complete and first countable.
Let $\{e_k\}_{k\in \mathbb I}$ be a Hilbert base for $\HC$.
\begin {enumerate}
\item If $U: \HC \to \KC$ is a unitary map then
      $\{f_k\}_{k\in \mathbb I}$ with $f_k := Ue_k$ is
      a Hilbert base for $\KC$.
\item If $\{f_k\}_{k\in \mathbb I}$ is a Hilbert base of $\KC$
      then there exists a unique unitary map $U: \HC \to \KC$
      such that $Ue_k = f_k$.
\end {enumerate}
\end {proposition}
For later use we finally mention the following corollary:
\begin {corollary} \label {GNSIdealRepCor}
Let $\HC$ be a pre-Hilbert space over
$\field C = \field R(i)$ where $\field R$ is real closed,
Cauchy complete and first countable and let $\{e_k\}_{k \in \mathbb I}$
be an orthonormal system such that
$\field C$-span$\{e_k\}_{k\in \mathbb I}$ is dense in $\HC$.
Denote the completion of $\HC$ by $\widehat{\HC}$. Then
$\{e_k\}_{k \in \mathbb I}$ is a Hilbert base of $\widehat{\HC}$ and
$\widehat{\HC}$ is unitary equivalent to $\ell^2(\field C)$, i.e.
there is a unitary map $U: \widehat{\HC} \to \ell^2(\field C)$.
\end {corollary}

We are now in the position to formulate a generalized GNS 
construction (analogous to the GNS pre-construction formulated in
proposition \ref {GNSConst}) taking into account the analytical 
properties of the field $\field C$. Using proposition 
\ref {CompleteHilbProp} to complete the GNS 
pre-Hilbert space (\ref{AquotientJ}) we get the following 
easy consequence of the GNS pre-construction 
(proposition \ref {GNSConst}):
\begin {theorem} [General GNS construction]
\label {FinalGNSConst}
Let $\field R$ be a real closed and Cau\-chy complete field and 
$\field C = \field R(i)$ and let
$\mathcal A$ be a $\field C$-algebra with involution $^*$.
Then for any positive linear functional 
$\omega: \mathcal A \to \field C$ there exists 
a Hilbert space $\hat{\mathfrak H}_\omega$ with a dense subspace
$\mathfrak H_\omega$ carrying a $^*$-representation
$\pi_\omega$ of $\mathcal A$ 
as constructed in proposition \ref {GNSConst}.
We shall call this representation the 
{\em GNS representation in $\hat{\mathfrak H}_\omega$}.
If $\mathcal A$ in addition has a unit element $\bfmath 1$ and
$\omega$ is a state then this representation is cyclic and we have
$\omega (A) = \SP{\psi_{\bfmath 1}}{\pi_\omega (A) \psi_{\bfmath 1}}$
and this property defines this representation up to 
unitary equivalence.
\end {theorem}

\section* {Acknowledgement}

The authors would like to thank K.~Fredenhagen for asking the question
of the GNS construction in deformation quantization and suggesting to
us to concentrate on the algebraic properties of $C^*$-algebras.
Moreover we would like to thank 
O.~Kegel for pointing out references \cite{GK77}, \cite{Gro79},
and \cite{Kell80}. Furthermore we would like to thank the referee 
for many useful and detailed comments and making us aware of 
references \cite {Han84,Kam86,Mai86,AV91,FO87,FW87}.
Finally, we would like to thank M.~Flato,
K.~Fredenhagen, J.~Huebschmann, H.~Rehren, and A.~Weinstein for
valuable discussions, and J.~Hoppe and M.~Lledo for a
charming table-football match.

\end {document}